\documentclass[aap,MSNbibl,seceqn,dvips]{arximspdf}

\doi{10.1214/13-AAP938} %kopijuoti is PTS
\volume{24}
\issue{3}
\pubyear{2014}
\firstpage{899}
\lastpage{934}

\makeatletter

\newcommand{\rright}{\right}
\newcommand{\lleft}{\left}
\newcommand{\rrVert}{\Vert}
\newcommand{\rrvert}{\vert}
\newcommand{\llVert}{\Vert}
\newcommand{\llvert}{\vert}
\newproclaim{st}{\textit{Step}}

\newtheorem{theorem}{Theorem}[section]
\newtheorem{proposition}[theorem]{Proposition}
\newtheorem{corollary}[theorem]{Corollary}
\newproclaim{assumption}[theorem]{Assumption}
\newtheorem{lemma}[theorem]{Lemma}
\newproclaim{remark}[theorem]{Remark}
\newproclaim{ass}{Assumption}
\makeatother

\begin{document}
\begin{frontmatter}

\title{Stochastic target games with controlled loss}
\runtitle{Stochastic target games with controlled loss}

\begin{aug}
\author[A]{\fnms{Bruno} \snm{Bouchard}\ead[label=e1]{bouchard@ceremade.dauphine.fr}\thanksref{t1}},
\author[B]{\fnms{Ludovic} \snm{Moreau}\ead[label=e2]{ludovic.moreau@math.ethz.ch}\thanksref{t2}}
\and
\author[C]{\fnms{Marcel} \snm{Nutz}\corref{}\ead[label=e3]{mnutz@math.columbia.edu}\thanksref{t3}}
\runauthor{B. Bouchard, L. Moreau and M. Nutz}
\affiliation{Universit\'e Paris Dauphine and CREST-ENSAE, ETH Z\"urich\\ and Columbia University}
\address[A]{B. Bouchard\\
CEREMADE and CREST-ENSAE\\
Universit\'e Paris Dauphine\\
75775 Paris Cedex 16\\
France\\
\printead{e1}}
\address[B]{L. Moreau\\
Department of Mathematics\\
ETH Z\"urich\\
8092 Z\"urich\\
Switzerland\\
\printead{e2}}
\address[C]{M. Nutz\\
Department of Mathematics\\
Columbia University\\
New York, New York 10027\\
USA\\
\printead{e3}}
\end{aug}
\thankstext{t1}{Supported by ANR Liquirisk and Labex Ecodec.}
\thankstext{t2}{Supported in part by European Research Council Grant 228053-FiRM and by the ETH Foundation.}
\thankstext{t3}{Supported by NSF Grant DMS-12-08985.}

% HISTORY:
\received{\smonth{6} \syear{2012}}
\revised{\smonth{4} \syear{2013}}

% ABSTRACT
%
\begin{abstract}
We study a stochastic game where one player tries to find a strategy
such that the state process reaches a target of controlled-loss-type,
no matter which action is chosen by the other player. We provide, in a
general setup, a relaxed geometric dynamic programming principle for
this problem and derive, for the case of a controlled SDE, the
corresponding dynamic programming equation in the sense of viscosity
solutions. As an example, we consider a problem of partial hedging
under Knightian uncertainty.
\end{abstract}

% KEYWORDS

\end{frontmatter}

\setcounter{footnote}{3}
%s1 #&#
\section{Introduction}

We study a stochastic (semi) game of the following form. Given an
initial condition $(t,z)$ in time and space, we try to find a strategy
$\mathfrak{u}[\cdot]$ such that the controlled state process
$Z_{t,z}^{\mathfrak{u}[\nu
],\nu}(\cdot)$ reaches a certain target at the given time $T$, no
matter which control $\nu$ is chosen by the adverse player. The target
is specified in terms of expected loss; that is, we are given a
real-valued (``loss'') function $\ell$ and try to keep the expected
loss above a given threshold $p\in\mathbb{R}$,
%
%
%e1.1 #&#
\begin{equation}
\label{eqtargetIntro} \mathop{\operatorname{ess}\inf}_\nu\mathbb{E} \bigl[
\ell \bigl(Z^{\mathfrak{u}[\nu],\nu
}_{t,z}(T) \bigr)|\mathcal{F} _{t}
\bigr]\geq p\qquad\mbox{a.s.}
\end{equation}
Instead of a game, one may also see this as a target problem under
Knightian uncertainty; then the adverse player has the role of choosing
a worst-case scenario.

Our aim is to describe, for given $t$, the set $\Lambda(t)$ of all
pairs $(z,p)$ such that there exists a strategy $\mathfrak{u}$
attaining the target.
We provide, in a general abstract framework, a geometric dynamic
programming principle~(GDP) for this set. To this end, $p$ is seen as
an additional state variable and formulated dynamically via a family
$\{M^\nu\}$ of auxiliary martingales with expectation $p$, indexed by
the adverse controls $\nu$. Heuristically, the GDP then takes the
following form: $\Lambda(t)$ consists of all $(z,p)$ such that there
exist a strategy $\mathfrak{u}$ and a family $\{M^\nu\}$ satisfying
\[
\bigl(Z^{\mathfrak{u}[\nu],\nu}_{t,z} (\tau ),M^{\nu
} (\tau ) \bigr)
\in{\Lambda} (\tau )\qquad\mbox{a.s.}
\]
for all adverse controls $\nu$ and all stopping times $\tau\geq t$. The
precise version of the GDP, stated in Theorem~\ref{thmGDP},
incorporates several relaxations that allow us to deal with various
technical problems. In particular, the selection of
$\varepsilon$-optimal strategies is solved by a covering argument which
is possible due to a continuity assumption on $\ell$ and a relaxation
in the variable $p$. The martingale $M^\nu$ is constructed from the
semimartingale decomposition of the adverse player's value process.

Our GDP is tailored such that the dynamic programming equation can be
derived in the viscosity sense. We exemplify this in
Theorem~\ref{thmpdederivavecz} for the standard setup where the state
process is determined by a stochastic differential equation (SDE) with
coefficients controlled by the two players; however, the general GDP
applies also in other situations such as singular control. The solution
of the equation, a partial differential equation (PDE) in our example,
corresponds to the \mbox{indicator} function of (the complement of) the
graph of $\Lambda$. In Theorem~\ref{thmpdederivavecx,y}, we specialize
to a case with a monotonicity condition, that is, particularly suitable
for pricing problems in mathematical finance. Finally, in order to
illustrate various points made throughout the paper, we consider a
concrete example of pricing an option with partial hedging, according
to a loss constraint, in a model where the drift and volatility
coefficients of the underlying are uncertain. In a worst-case analysis,
the uncertainty corresponds to an adverse player choosing the
coefficients; a formula for the corresponding seller's price is given
in Theorem~\ref{thmexplicitexpression}.

Stochastic target (control) problems with almost-sure constraints,
corresponding to the case where $\ell$ is an indicator function and
$\nu $ is absent, were introduced in \cite{SonerTouzi02b,SonerTouzi02a}
as an extension of the classical superhedging
problem~\cite{ElKarouiQuenez95} in mathematical finance. Stochastic
target problems with controlled loss were first studied in
\cite{BouchardElieTouzi09} and are inspired by the quantile hedging
problem~\cite{FollmerLeukert99}. The present paper is the first to
consider stochastic target games.
The rigorous treatment of zero-sum stochastic differential games was
pioneered in \cite{FlemingSouganidis89}, where the mentioned selection
problem for $\varepsilon$-optimal strategies was treated by a
discretization and a passage to continuous-time limit in the PDEs. Let
us remark, however, that we have not been able to achieve satisfactory
results for our problem using such techniques. We have been importantly
influenced by \cite{BuckdahnLi08}, where the value functions are
defined in terms of essential infima and suprema, and then shown to be
deterministic. The formulation with an essential infimum (rather than
an infimum of suitable expectations) in~(\ref{eqtargetIntro}) is
crucial in our case, mainly because $\{M^\nu\}$ is constructed by a
method of non-Markovian control, which raises the fairly delicate
problem of dealing with one nullset for every adverse control $\nu$.

The remainder of the paper is organized as follows.
Section~\ref{sectgdp} contains the abstract setup and GDP. In
Section~\ref{secexample} we specialize to the case of a controlled SDE
and derive the corresponding PDE, first in the general case and then in
the monotone case. The problem of hedging under uncertainty is
discussed in Section~\ref{sectpartialHedging}.

%s2 #&#
\section{Geometric dynamic programming principle}\label{sectgdp}

In this section, we obtain our geometric dynamic programming principle
(GDP) in an abstract framework. Some of our assumptions are simply the
conditions we need in the proof of the theorem; we will illustrate
later how to actually verify them in a typical setup.

%s2.1 #&#
\subsection{Problem statement}

We fix a time horizon $T>0$ and a probability space $(\Omega,\mathcal
{F},\mathbb{P})$ equipped with a filtration
$\mathbb{F}=(\mathcal{F}_{t})_{t\in [0,T]}$ satisfying the usual
conditions of right-continuity and completeness. We shall consider two
sets $\mathcal{U}$ and $\mathcal{V}$ of controls; for the sake of
concreteness, we assume that each of these sets consists of stochastic
processes on $(\Omega,\mathcal{F})$, indexed by $[0,T]$ and with values
in some sets $U$~and~$V$, respectively. Moreover, let $\mathfrak{U}$ be
a set of mappings $\mathfrak{u}\dvtx  \mathcal{V}\to\mathcal{U}$. Each
$\mathfrak {u}\in\mathfrak{U}$ is called a strategy, and the notation
$\mathfrak{u}[\nu]$ will be used for the control it associates with
$\nu\in\mathcal{V}$. In applications, $\mathfrak{U}$ will be chosen to
consist of mappings that are nonanticipating; see
Section~\ref{secexample} for an example. Furthermore, we are given a
metric space $(\mathcal{Z},d_{\mathcal {Z}})$ and, for each
$(t,z)\in[0,T]\times\mathcal{Z}$ and $(\mathfrak{u},\nu) \in
\mathfrak{U} \times\mathcal{V}$, an adapted c\` adl\`ag process
$Z^{\mathfrak{u}[\nu],\nu}_{t,z}(\cdot)$ with values in $\mathcal{Z}$
satisfying $Z_{t,z}^{\mathfrak{u}[\nu],\nu}(t)=z$. For brevity, we set
\[
Z^{\mathfrak{u},\nu}_{t,z}:=Z^{\mathfrak{u}[\nu],\nu}_{t,z}.
\]
Let $\ell\dvtx \mathcal{Z}\to\mathbb{R}$ be a Borel-measurable function
satisfying
%
%e2.1 #&#
\begin{equation}\label{eqelldansL1}
\mathbb{E} \bigl[\bigl\llvert \ell \bigl(Z^{\mathfrak{u},\nu
}_{t,z}(T)
\bigr)\bigr\rrvert \bigr]<\infty\qquad\mbox{for all }(t,z,\mathfrak{u},\nu)
\in[0,T]\times\mathcal {Z}\times\mathfrak{U}\times\mathcal{V}.
\end{equation}
We interpret
$\ell$ as a loss (or ``utility'') function and denote by
\[
I(t,z,\mathfrak{u},\nu):=\mathbb{E} \bigl[\ell \bigl(Z^{\mathfrak
{u},\nu}_{t,z}(T)
\bigr)|\mathcal{F} _{t} \bigr]\qquad(t,z,\mathfrak{u},\nu)\in[0,T]
\times\mathcal {Z}\times\mathfrak{U}\times\mathcal{V}
\]
the expected loss given $\nu$ (for the player choosing $\mathfrak
{u}$) and by
\[
J(t,z,\mathfrak{u}):=\mathop{\operatorname{ess}\inf}_{\nu\in
\mathcal{V}} I(t,z,
\mathfrak {u},\nu)\qquad(t,z,\mathfrak{u})\in [0,T]\times\mathcal{Z}\times
\mathfrak{U}
\]
the worst-case expected loss.
The main object of this paper is the reachability set
%
%e2.2 #&#
\begin{eqnarray}\label{eqdefinitiondureachabilityset}
\Lambda(t)&:=& \bigl\{(z,p) \in\mathcal{Z}\times\mathbb{R}: \mbox{ there exists }\mathfrak{u}\in\mathfrak{U}
\nonumber\\[-8pt]\\[-8pt]
&&\hspace*{26pt} \mbox{ such that }J(t,z,\mathfrak{u})
\ge p\mbox{ $\mathbb{P}$-a.s.} \bigr\}.\nonumber
\end{eqnarray}
These are the initial conditions $(z,p)$ such that starting at time
$t$, the player choosing $\mathfrak{u}$ can attain an expected loss not
worse than $p$, regardless of the adverse player's action $\nu$. The
main aim of this paper is to provide a geometric dynamic programming
principle for $\Lambda(t)$. For the case without adverse player, a
corresponding result was obtained in~\cite{SonerTouzi02b} for the
target problem with almost-sure constraints and
in~\cite{BouchardElieTouzi09} for the problem with controlled loss.

As mentioned above, the dynamic programming for
problem~(\ref{eqdefinitiondureachabilityset}) requires the introduction
of a suitable set of martingales starting from $p\in\mathbb{R}$. This
role will be played by certain families\footnote{Of course, there is no
mathematical difference between families indexed by $\mathcal{V}$, like
$\{M^{\nu},\nu\in\mathcal{V}\}$, and mappings on $\mathcal{V}$, like
$\mathfrak{u}$. We shall use both notions interchangeably, depending on
notational convenience.} $\{M^{\nu},\nu\in\mathcal{V}\}$ of martingales
which should be considered as additional controls. More precisely, we
denote by $\mathcal{M}_{t,p}$ the set of all real-valued
(right-continuous) martingales $M$ satisfying $M(t)=p$
$\mathbb{P}$-a.s., and we\vadjust{\goodbreak} fix a set
$\mathfrak{M}_{t,p}$ of families $\{M^{\nu
},\nu\in\mathcal{V}\}\subset\mathcal{M}_{t,p}$; further assumptions on
$\mathfrak{M}_{t,p}$ will be introduced below. Since these martingales
are not present in the original
problem~(\ref{eqdefinitiondureachabilityset}), we can choose
$\mathfrak{M}_{t,p}$ at our convenience; see also
Remark~\ref{remmartIsNonanticip} below.

As usual in optimal control, we shall need to concatenate controls and
strategies in time according to certain events. We use the notation
\[
\nu\oplus_{\tau} \bar{\nu}:= \nu
\mathbf{1}_{[0,\tau]}+\bar {\nu}\mathbf {1}_{(\tau,T]}
\]
for the concatenation of two controls $\nu,\bar{\nu}\in\mathcal
{V}$ at a
stopping time $\tau$. We also introduce the set
\[
\{\nu=_{{(t,\tau]}}\bar{\nu} \}:= \bigl\{\omega\in\Omega\dvtx  \nu
_s(\omega)=\bar{\nu}_s(\omega)\mbox{ for all }s\in\bigl(t,
\tau(\omega )\bigr] \bigr\}.
\]
Analogous notation is used for elements of $\mathcal{U}$.

In contrast to the setting of control, strategies can be concatenated
only at particular events and stopping times, as otherwise the
resulting strategies would fail to be elements of $\mathfrak{U}$ (in
particular, because they may fail to be nonanticipating, see also
Section~\ref{secexample}). Therefore, we need to formalize the events
and stopping times which are admissible for this purpose: for each
$t\le T$, we consider a set $\mathfrak{F}_{t}$ whose elements are
families $\{A^{\nu },\nu\in\mathcal{V}\}\subset\mathcal{F}_{t}$ of
events indexed by $\mathcal{V}$, as well as a set $\mathfrak{T}_{t}$
whose elements are families $\{\tau^{\nu},\nu \in \mathcal{V}\}
\subset\mathcal{T}_{t}$, where $\mathcal{T}_{t}$ denotes the set of all
stopping times with values in $[t,T]$. We assume that
$\mathfrak{T}_{t}$ contains any deterministic time $s\in[t,T]$ (seen as
a constant family $\tau^\nu \equiv s$, $\nu\in\mathcal{V}$). In\vadjust{\goodbreak}
practice, the sets $\mathfrak{F}_{t}$ and $\mathfrak{T}_{t}$ will not
contain all families of events and stopping times, respectively; one
will impose additional conditions on $\nu\mapsto A^\nu$ and
$\nu\mapsto\tau ^\nu$ that are compatible with the conditions defining
$\mathfrak{U}$. Both sets should be seen as auxiliary objects which
make it easier (if not possible) to verify the dynamic programming
conditions below.

%
%s2.2 #&#
\subsection{The geometric dynamic programming principle}

We can now state the conditions for our main result. The first one
concerns the concatenation of controls and strategies.

\renewcommand{\theass}{(C)}
%
%as1 #&#
\begin{ass}\label{assC}
The following hold for all $ t \in[0,T]$:
\begin{longlist}[(C6)]
\item[(C1)]\label{assVcconcat} Fix $\nu_{0},\nu_{1},\nu_{2}\in
    \mathcal{V}$ and $A\in\mathcal{F}_{t} $. Then $\nu:= \nu_{0}
    \oplus_{t} ( \nu _{1}{\mathbf
    1}_{A}+\nu_{2}{\mathbf1}_{A^{c}})\in\mathcal{V}$.

\item[(C2)]\label{assUcconcat} Fix $(\mathfrak{u}_{j})_{j\ge0}\subset
    \mathfrak{U}$, and let $\{ A^{\nu}_{j},\nu\in\mathcal{V}\}_{j\ge1}
    \subset \mathfrak{F} _{t} $ be such that $\{A^{\nu}_{j},j\ge1\} $
    forms a partition of $\Omega$ for each $\nu\in\mathcal{V}$. Then
    $\mathfrak{u}\in\mathfrak{U}$ for
\[
\mathfrak{u}[\nu]:= \mathfrak{u}_{0}[\nu] \oplus
_{t} \sum_{j\ge
1}
\mathfrak{u}_{j}[\nu] \mathbf{1}_{A_{j}^{\nu}},\qquad\nu\in
\mathcal{V}.
\]

\item[(C3)]\label{assUcWithFixedNu} Let $\mathfrak{u}\in\mathfrak{U}$
    and
    $\nu \in\mathcal{V}$. Then
    $\mathfrak{u}[\nu\oplus_t\cdot]\in\mathfrak{U}$.

\item[(C4)]\label{assAdmiA} Let $\{A^{\nu},\nu\in\mathcal{V}\}
    \subset\mathcal{F}_{t}$ be a family of events such that
    $A^{\nu_{1}}\cap\{\nu_{1}=_{{(0,t]}} \nu_{2}\}=
    A^{\nu_{2}}\cap\{\nu_{1}=_{{ (0,t]}} \nu_{2}\}$ for all $\nu
    _{1},\nu_{2}\in\mathcal{V}$. Then $\{A^{\nu},\nu\in\mathcal{V}\}
    \in\mathfrak{F}_{t}$.

\item[(C5)]\label{asstauegauxsicontroleegauxsurevenement} Let
    $\{\tau^{\nu},\nu\in\mathcal{V}\}\in\mathfrak{T}_{t}$. Then
    $\{\tau ^{\nu_{1}}\le s\}\cap\{\nu_{1}=_{(0,s]}\nu_{2}\}=
    \{\tau^{\nu_{2}}\le s\}\cap\{ \nu _{1}=_{(0,s]}\nu_{2}\}$
    $\mathbb{P}$-a.s. for all $\nu_{1},\nu _{2}\in\mathcal{V}$ and
    $s\in[t,T]$.

\item[(C6)]\label{assAdmitau} Let $\{\tau^{\nu},\nu\in\mathcal{V}\}
    \in\mathfrak{T}_{t}$. Then, for all $t\le s_{1}\le s_{2}\le T$,
    $\{\{ \tau^{\nu}\in (s_{1},s_{2}]\},\nu\in\mathcal{V}\}$ and
    $\{\{\tau^{\nu}\notin (s_{1},s_{2}]\},\nu\in\mathcal{V}\}$ belong
    to $\mathfrak{F}_{s_{2}}$.\vadjust{\goodbreak}
\end{longlist}
\end{ass}

The second condition concerns the behavior of the state process.

\renewcommand{\theass}{(Z)}
%
%as2 #&#
\begin{ass}\label{assZ}
The following hold for all $(t,z,p)\in[0,T]\times\mathcal{Z}\times
\mathbb{R}$ and $s\in[t,T]$:
\begin{longlist}[(Z4)]
\item[(Z1)]\label{assZegauxsistrategalessurevenement} $
    Z^{\mathfrak{u}_{1},\nu}_{t,z}(s)(\omega)=Z^{\mathfrak{u}_{2},\nu
    }_{t,z}(s)(\omega)$ for $\mathbb{P}$-a.e.
    $\omega\in\{\mathfrak{u}_{1}[\nu
    ]=_{(t,s]}\mathfrak{u}_{2}[\nu]\}$, for all $\nu\in\mathcal{V}$ and
    $\mathfrak{u}_{1},\mathfrak{u}_{2}\in \mathfrak{U}$.

\item[(Z2)]\label{assZegauxsicontroleegauxsurevenement} $
    Z^{\mathfrak{u},\nu_{1}}_{t,z}(s)(\omega)=Z^{\mathfrak{u},\nu
    _{2}}_{t,z}(s)(\omega)$ for $\mathbb{P}$-a.e.
    $\omega\in\{\nu_{1}=_{(0,s]}\nu_{2}\}$, for all $\mathfrak{u}
    \in\mathfrak{U}$ and $\nu_{1},\nu_{2}\in\mathcal{V}$.

\item[(Z3)]\label{assMegauxsicontroleegauxsurevenement} $
    M^{\nu_{1}}(s)(\omega)=M^{\nu_{2}}(s)(\omega)$ for $\mathbb
    {P}$-a.e. $\omega \in\{\nu_{1}=_{(0,s]}\nu_{2}\}$, for all
    $\{M^{\nu},\nu\in \mathcal{V}\}\in \mathfrak{M}_{t,p}$ and
    $\nu_{1},\nu_{2}\in\mathcal{V}$.

\item[(Z4)]\label{assesssupessinf=constant} There exists a constant
    $K(t,z)\in\mathbb{R}$ such that
\[
\mathop{\operatorname{ess}\sup}_{\mathfrak{u}\in\mathfrak{U}} \mathop{\operatorname{ess}
\inf}_{ \nu\in\mathcal{V}} \mathbb{E} \bigl[\ell\bigl(Z^{\mathfrak{u}, \nu
}_{t,z}(T)
\bigr)|\mathcal{F}_{t} \bigr]=K(t,z)\qquad\mathbb{P}\mbox{-a.s.}
\]
\end{longlist}
\end{ass}

The nontrivial assumption here is, of course,
\hyperref[assesssupessinf=constant]{\textup{(Z4)}}, stating that (a
version of) the random variable
$\mathop{\operatorname{ess}\sup}_{\mathfrak{u}\in\mathfrak{U}}
\mathop{\operatorname{ess}\inf}_{ \nu\in\mathcal{V}}
\mathbb{E}[\ell(Z^{\mathfrak{u}, \nu }_{t,z}(T))|\mathcal{F}_{t}]$ is
\textit{deterministic}. For the game determined by a Brownian SDE as
considered in Section~\ref{secexample}, this will be true by a result
of~\cite{BuckdahnLi08}, which, in turn, goes back to an idea of
\cite{Peng97b} (see also \cite{LiPeng09}). An extension to jump
diffusions can be found in \cite{BuckdahnHuLi11}.

While the above assumptions are fundamental, the following conditions
are of technical nature. We shall illustrate later how they can be verified.

\renewcommand{\theass}{(I)}
%
%as3 #&#
\begin{ass}\label{assI} Let $(t,z)\in[0,T]\times\mathcal{Z}$,
$\mathfrak{u}\in\mathfrak{U}$ and $\nu \in\mathcal{V}$.
\begin{longlist}[(I2)]
\item[(I1)]\label{assIborneinf} There exists an adapted
    right-continuous
    process $N^{\mathfrak{u},\nu}_{t,z}$ of class~(D) such that
\[
\mathop{\operatorname{ess}\inf}_{\bar\nu\in\mathcal{V}} \mathbb{E} \bigl[\ell
\bigl(Z_{t,z}^{\mathfrak{u},\nu\oplus_{s}
\bar\nu}(T) \bigr)|\mathcal{F}_{s}
\bigr] \ge N^{\mathfrak{u},\nu
}_{t,z}(s) \qquad\mbox{$\mathbb{P}$-a.s. for all $s\in[t,T]$}.
\]

\item[(I2)]\label{assIbornesup} There exists an adapted
    right-continuous
    process $L^{\mathfrak{u},\nu}_{t,z}$ such that $L^{\mathfrak{u},\nu
    }_{t,z}(s)\in L^{1}$ and
\[
\mathop{\operatorname{ess}\inf}_{\bar\mathfrak{u}\in\mathfrak
{U}} \mathbb{E} \bigl[\ell
\bigl(Z_{t,z}^{\mathfrak{u}\oplus_{s}\bar\mathfrak{u},\nu}(T) \bigr)|\mathcal{F}_{s}
\bigr] \ge L^{\mathfrak{u},\nu}_{t,z}(s)\qquad\mbox{$\mathbb{P}$-a.s. for all $s\in[t,T]$}.
\]
Moreover, $ L^{\mathfrak{u},\nu_{1}}_{t,z}(s)(\omega)=L^{\mathfrak
{u},\nu _{2}}_{t,z}(s)(\omega)$ for $\mathbb{P}$-a.e. $\omega\in\{\nu
_{1}=_{(0,s]}\nu _{2}\}$, for all $\mathfrak{u}\in\mathcal{U}$ and
$\nu_{1},\nu _{2}\in\mathcal{V}$.
\end{longlist}
\end{ass}

\renewcommand{\theass}{(R)}
%
%as4 #&#
\begin{ass}\label{assR}
Let $(t,z)\in[0,T]\times\mathcal{Z}$.
\begin{longlist}[(R2)]
\item[(R1)]\label{assreguIJK} Fix $s\in[t,T]$ and $\varepsilon>0$. Then
    there exist a Borel-measurable partition $(B_{j})_{j\geq1}$ of $\mathcal {Z}$ and
    a sequence $(z_{j})_{j\ge1}\subset\mathcal{Z}$ such that for all
    $\mathfrak{u}\in\mathfrak{U}$, $\nu\in\mathcal{V}$ and~$j\ge1$,
%
%e2.3 #&#
\begin{eqnarray}
\left.
\begin{array}{rcl}
\mathbb{E} \bigl[\ell\bigl(Z^{\mathfrak{u},\nu}_{t,z}(T)
\bigr)|\mathcal {F}_{s} \bigr]& \ge& I(s,z_{j},\mathfrak{u},
\nu) - \varepsilon,
\vspace*{2pt}\cr
\displaystyle\mathop{\operatorname{ess}\inf}_{\bar\nu\in\mathcal{V}} \mathbb {E} \bigl[\ell
\bigl(Z^{\mathfrak{u},\nu\oplus_{s}\bar\nu
}_{t,z}(T)\bigr)|\mathcal{F}_{s}
\bigr] &\le& J\bigl(s,z_{j},\mathfrak {u}[\nu\oplus
_s\cdot]\bigr) + \varepsilon,
\cr
K(s,z_{j})-\varepsilon &\le& K\bigl(s,Z^{\mathfrak{u},\nu}_{t,z}(s)
\bigr) \le K(s,z_{j})+\varepsilon
\end{array}
\right\}\nonumber
\\
\eqntext{\mathbb{P}\mbox{-a.s. on
}\bigl\{Z^{\mathfrak{u},\nu}_{t,z}(s)\in B_{j}\bigr\}.}
\end{eqnarray}

\item[(R2)]\label{assreguprobaZ} $
    \lim_{\delta\to0}\sup_{\nu\in\mathcal{V}, \tau\in\mathcal {T}_{t}}
    \mathbb{P} \{{\sup_{0\le h\le\delta}{ d_{\mathcal{Z}}}
    (Z^{\mathfrak{u}, \nu
    }_{t,z}(\tau+h),Z^{\mathfrak{u},\nu}_{t,z}(\tau) )\ge\varepsilon}
    \} = 0$ for all $\mathfrak{u}\in\mathfrak{U}$ and $\varepsilon>0$.
\end{longlist}
\end{ass}

Our GDP will be stated in terms of the closure
\[
\bar\Lambda(t):= \lleft\{\begin{array}{l}
(z,p)\in\mathcal{Z}\times\mathbb{R}:
\mbox{ there exist }(t_{n},z_{n},p_{n})\to(t,z,p)
\\
\mbox{such that } (z_{n},p_{n})\in\Lambda(t_{n})
\mbox{ and }t_{n}\ge t\mbox{ for all } n\ge1
\end{array} \rright\}
\]
and the uniform interior
\[
\Lambda_{\iota}(t):= \bigl\{(z,p)\in\mathcal {Z}\times
\mathbb{R}\dvtx  \bigl(t',z',p'\bigr)\in
B_{\iota}(t,z,p)\mbox{ implies } \bigl(z',p'
\bigr)\in\Lambda\bigl(t'\bigr) \bigr\},
\]
where $B_{\iota}(t,z,p)\subset[0,T]\times\mathcal{Z}\times \mathbb{R}$
denotes the open ball with center $(t,z,p)$ and radius $\iota>0$ [with
respect to the distance function
$d_{\mathcal{Z}}(z,z')+|p-p'|+|t-t'|$]. The relaxation from $\Lambda$
to $\bar\Lambda$ and $\Lambda _{\iota}$ essentially allows us to reduce
to stopping times with countably many values in the proof of the GDP
and thus to avoid regularity assumptions in the time variable. We shall
also relax the variable $p$ in the assertion of (GDP2); this is
inspired by~\cite{BouchardNutz11} and important for the covering
argument in the proof of (GDP2), which, in turn, is crucial due to the
lack of a measurable selection theorem for strategies. Of course, all
our relaxations are tailored such that they will not interfere
substantially with the derivation of the dynamic programming equation;
cf. Section~\ref{secexample}.

%
%th2.1 #&#
\begin{theorem}\label{thmGDP}
Fix $(t,z,p)\in[0,T]\times\mathcal{Z}\times\mathbb{R}$ and let
Assumptions \ref{assC}, \ref{assZ}, \ref{assI} and \ref{assR} hold
true.
\begin{longlist}[(GDP1)]
\item[(GDP1)] If $(z,p)\in\Lambda(t)$, then there exist
    $\mathfrak{u}\in \mathfrak{U}$ and $\{M^{\nu},\nu\in{\mathcal{V}}
    \}\subset \mathcal {M}_{t,p}$ such that
\[
\bigl(Z^{\mathfrak{u},\nu}_{t,z} (\tau ),M^{\nu} (\tau ) \bigr)
\in \bar\Lambda (\tau ) \qquad\mathbb{P}\mbox{-a.s. for all $\nu\in
\mathcal{V}$ and $\tau\in\mathcal{T}_{t}$}.
\]

\item[(GDP2)] Let $\iota>0$, $\mathfrak{u}\in\mathfrak
    {U}$, $\{ M^{\nu},\nu\in\mathcal{V}\} \in\mathfrak{M}_{t,p}$ and
    $\{\tau^{ \nu}, \nu\in\mathcal{V}\} \in\mathfrak{T} _{t}$ be such
    that
\[
\bigl(Z^{\mathfrak{u},\nu}_{t,z}\bigl(\tau^{\nu}
\bigr),M^{\nu}\bigl(\tau^{\nu
}\bigr) \bigr) \in
\Lambda_{\iota}\bigl(\tau^{\nu}\bigr) \qquad\mathbb{P}\mbox{-a.s. for all } \nu \in\mathcal{V}
\]
and suppose that $ \{M^{\nu}(\tau^\nu)^{+}\dvtx  \nu\in\mathcal {V} \}$
and $ \{L^{\mathfrak{u},\nu}_{t,z}(\tau')^{-}\dvtx  \nu\in\mathcal {V},
\tau'\in\mathcal{T}_t \}$ are uniformly integrable, where
$L^{\mathfrak{u},\nu}_{t,z}$ is as
in~\hyperref[assIbornesup]{\textup{(I2)}}. Then
$(z,p-\varepsilon)\in\Lambda(t)$ for all $\varepsilon>0$.
\end{longlist}
\end{theorem}

The proof is stated in Sections~\ref{secproofGDP1}
and~\ref{secproofGDP2} below.

%
%re2.2 #&#
\begin{remark}\label{remmartIsNonanticip}
We shall see in the proof that the family $\{M^{\nu},\nu\in{\mathcal
{V}} \} \subset\mathcal{M}_{t,p}$ in (GDP1) can actually be chosen to
be nonanticipating in the sense
of~\hyperref[assMegauxsicontroleegauxsurevenement]{\textup{(Z3)}}.
However, this will not be used when (GDP1) is applied to derive the
dynamic programming equation. Whether $\{M^{\nu},\nu\in{\mathcal{V}}
\}$ is an element of $\mathfrak{M}_{t,p}$ will depend on the definition
of the latter set; in fact, we did not make any assumption about its
richness. In many application, it is possible to take
$\mathfrak{M}_{t,p}$ to be the set of all nonanticipating families in
$\mathcal{M}_{t,p}$; however, we prefer to leave some freedom for the
definition of $\mathfrak{M}_{t,p}$ since this may be useful in ensuring
the uniform integrability required in~(GDP2).
\end{remark}

We conclude this section with a version of the GDP for the case
$\mathcal{Z}=\mathbb{R} ^{d}$, where we show how to reduce from
standard regularity conditions on the state process and the loss
function to the Assumptions
\hyperref[assreguIJK]{\textup{(R1)}}~and~\ref{assI}.

%
%co2.3 #&#
\begin{corollary} \label{corGDPwithlessregularity}
Let Assumptions \ref{assC}, \ref{assZ} and
\hyperref[assreguprobaZ]{\textup{(R2)}} hold true. Assume also that
$\ell$ is continuous and that there exist constants $C\geq0$ and
$\bar{q}>q\geq0$ and a~locally bounded function $\varrho\dvtx
\mathbb{R}^d\mapsto\mathbb{R}_{+}$ such that
%
%
%e2.4 #&#
%e2.5 #&#
\begin{eqnarray}
\bigl|\ell(z)\bigr| &\le& C\bigl(1+|z|^{q}\bigr), \label{eqassellZpoly}
\\
\mathop{\operatorname{ess}\sup}_{(\bar\mathfrak{u},\bar\nu)\in
\mathfrak{U}\times\mathcal{V}} \mathbb{E}
\bigl[\bigl|Z_{t,z}^{ \bar\mathfrak{u},
\bar\nu}(T)\bigr|^{\bar{q}}|\mathcal{F}_{t}
\bigr] &\le&\varrho( z)^{\bar{q}} \qquad\mathbb{P}\mbox{-a.s.}\label{eqassZgrowth}
\end{eqnarray}
and
%e2.6 #&#
\begin{eqnarray}\label{eqZlipschitz}
&& \mathop{\operatorname{ess}\sup}_{(\bar\mathfrak{u},\bar\nu)\in
\mathfrak{U}\times\mathcal{V}} \mathbb{E}
\bigl[\bigl|Z_{t,z}^{\mathfrak{u}\oplus_{s}\bar\mathfrak{u},\nu\oplus_{s}\bar\nu}(T)-Z_{s,z'}^{\bar
\mathfrak{u},\nu\oplus_s\bar
\nu}(T)\bigr| |
\mathcal{F}_{s} \bigr]
\nonumber\\[-8pt]\\[-8pt]
&&\qquad \le C\bigl|Z_{t,z}^{\mathfrak{u},\nu
}(s)-z'\bigr|\qquad\mathbb{P}\mbox{-a.s.}\nonumber
\end{eqnarray}
for all $(t,z)\in[0,T]\times\mathbb{R}^d$, $(s,z')\in[t,T]\times
\mathbb{R}^{d}$ and $(\mathfrak{u},\nu
)\in\mathfrak{U}\times\mathcal{V}$.

Let $(t,z)\in[0,T]\times\mathbb{R}^d$, and let $\{\tau^{\mathfrak
{u},\nu}, (\mathfrak{u},\nu)\in\mathfrak{U} \times\mathcal{V}\}
\subset\mathcal{T}_{t}$ be such that the collection
$\{Z^{\mathfrak{u},\nu }_{t,z} (\tau^{\mathfrak{u},\nu}
),(\mathfrak{u},\nu )\in\mathfrak{U}\times\mathcal{V}\}$ is uniformly
bounded in $L^{\infty}$.
\begin{longlist}[(GDP1$'$)]
\item[(GDP1$'$)] If $(z,p+\varepsilon)\in\Lambda(t)$ for
    some $\varepsilon>0$, then there exist
    $\mathfrak{u}\in\mathfrak{U}$ and $\{M^{\nu},\nu \in {\mathcal{V}}
    \}\subset \mathcal{M}_{t,p}$ such that
\[
\bigl(Z^{\mathfrak{u},\nu}_{t,z}\bigl(\tau^{\mathfrak{u},\nu}
\bigr),M^{\nu
}\bigl(\tau^{\mathfrak{u},\nu}\bigr) \bigr) \in \bar\Lambda
\bigl(\tau^{\mathfrak{u},\nu}\bigr) \qquad\mathbb{P}\mbox {-a.s. for all $\nu\in
\mathcal{V}$.}
\]

\item[(GDP2$'$)] If $\iota>0$, $\mathfrak{u}\in\mathfrak{U}$ and $\{
M^{\nu},\nu\in \mathcal{V}\}\in\mathfrak{M}_{t,p}$ are such that
\[
\bigl(Z^{\mathfrak{u},\nu}_{t,z}\bigl(\tau^{\mathfrak{u},\nu}
\bigr),M^{\nu
}\bigl(\tau^{\mathfrak{u},\nu}\bigr) \bigr) \in \Lambda_{\iota}\bigl(\tau^{\mathfrak{u},\nu}\bigr) \qquad\mathbb{P}
\mbox{-a.s. for all } \nu\in\mathcal{V}
\]
and $\{\tau^{\mathfrak{u},\nu}, \nu\in\mathcal{V}\} \in\mathfrak{T}
_{t}$, then $(z,p-\varepsilon)\in
\Lambda(t)$ for all $\varepsilon>0$.
\end{longlist}
\end{corollary}

We remark that Corollary~\ref{corGDPwithlessregularity} is usually
applied in a setting where $\tau^{\mathfrak{u}, \nu}$ is the exit time
of $Z^{\mathfrak{u},\nu}_{t,z}$ from a given ball, so that the
boundedness assumption is not restrictive. (Some adjustments are needed
when the state process admits unbounded jumps; see also
\cite{Moreau10}.)

%s2.3 #&#
\subsection{Proof of (GDP1)}\label{secproofGDP1}

We fix $t\in[0,T]$ and $(z,p)\in\Lambda(t)$ for the remainder of this
proof. By definition~(\ref{eqdefinitiondureachabilityset}) of
$\Lambda(t)$, there exists $\mathfrak{u}\in\mathfrak{U}$ such that
%
%e2.7 #&#
\begin{eqnarray}
\label{eqproofGDP1atteintcible} \qquad\mathbb{E} \bigl[G(\nu)|\mathcal{F}_{t} \bigr]&\ge& p
\qquad\mathbb {P}\mbox{-a.s. for all }\nu\in\mathcal{V},\mbox{ where
}G(\nu):=\ell\bigl(Z_{t,z}^{\mathfrak{u},\nu}(T)\bigr).
\end{eqnarray}
In order to construct the family $\{M^{\nu},\nu\in{\mathcal{V}} \}
\subset\mathcal{M}
_{t,p}$ of martingales, we consider
%
%e2.8 #&#
\begin{eqnarray}
\label{eqdefSnu} S^{\nu}(r)&:=&\mathop{\operatorname{ess}
\inf}_{\bar\nu\in\mathcal
{V}} \mathbb{E} \bigl[G(\nu\oplus
_{r}\bar\nu )|\mathcal{F}_{r} \bigr],
\qquad t\le r\le T.
\end{eqnarray}
We shall obtain $M^{\nu}$ from a Doob--Meyer-type decomposition of
$S^\nu$. This can be seen as a generalization with respect
to~\cite{BouchardElieTouzi09}, where the necessary martingale was
trivially constructed by taking the conditional expectation of the
terminal reward.%\vspace*{6pt}

\textit{Step} 1. We have $S^\nu(r) \in L^1(\mathbb{P})$ and $\mathbb{E}
[S^{\nu}(r)|\mathcal{F}_{s} ]\ge S^{\nu}(s)$ for all $t\leq s\leq r\leq
T$ and $\nu\in\mathcal {V}$.%\vspace*{6pt}

The integrability of $S^\nu(r)$ follows from (\ref{eqelldansL1}) and
\hyperref[assIborneinf]{\textup{(I1)}}. To see the submartingale
property, we first show that the family
$\{\mathbb{E}[G(\nu\oplus_{r}\bar\nu )|\mathcal{F}_{r}],
\bar\nu\in\mathcal{V}\}$ is directed downward. Indeed, given
$\bar\nu_{1},\bar\nu_{2}\in\mathcal{V}$, the set
\[
A:= \bigl\{\mathbb{E} \bigl[G(\nu\oplus_{r}
\bar\nu_1)|\mathcal {F}_{r} \bigr] \le\mathbb{E}
\bigl[G(\nu\oplus_{r}\bar
\nu_2)|\mathcal{F}_{r} \bigr] \bigr\}
\]
is in $\mathcal{F}_{r}$; therefore, $\bar\nu_{3}:=\nu\oplus_{r}
(\bar\nu_{1} \mathbf{1}_{A}+\bar\nu_{2} \mathbf{1}_{A^c})$ is an
element of $\mathcal{V}$ by
Assumption~\hyperref[assVcconcat]{\textup{(C1)}}. Hence,
\hyperref[assZegauxsicontroleegauxsurevenement]{\textup{(Z2)}} yields
that
\begin{eqnarray*}
\mathbb{E} \bigl[G(\nu\oplus_{r}\bar
\nu_3)|\mathcal{F}_{r} \bigr] &=&\mathbb{E}
\bigl[G(\nu\oplus_{r}\bar
\nu_{1})\mathbf {1}_{A}+G(\nu\oplus
_{r}\bar \nu_{2})\mathbf{1}_{A^c}
|\mathcal{F}_{r} \bigr]
\\
&=& \mathbb{E} \bigl[G(\nu\oplus_{r}
\bar\nu_1)|\mathcal {F}_{r} \bigr]
\mathbf{1}_{A}+\mathbb{E} \bigl[G(\nu \oplus
_{r}\bar\nu_2)|\mathcal{F}_{r} \bigr] \mathbf{1}_{A^c}
\\
&=& \mathbb{E} \bigl[G(\nu\oplus_{r}
\bar\nu_1)|\mathcal {F}_{r} \bigr] \wedge\mathbb{E}
\bigl[G(\nu\oplus_{r}\bar
\nu_2)|\mathcal{F}_{r} \bigr].
\end{eqnarray*}
As a result, we can find a sequence $(\bar\nu_{n})_{n\ge1}$ in
$\mathcal{V}$ such that
$\mathbb{E}[G(\nu\oplus_{r}\bar\nu_n)|\mathcal{F}_{r}]$ decreases
$\mathbb{P}$-a.s. to $S^\nu(r)$; cf. \cite{Neveu75},
Proposition~VI-1-1. Recalling (\ref{eqelldansL1}) and that $S^\nu(r)\in
L^1(\mathbb{P})$, monotone convergence yields that
\begin{eqnarray*}
\mathbb{E} \bigl[S^{\nu}(r)|\mathcal{F}_{s} \bigr] &=&
\mathbb{E} \Bigl[\lim_{n\to\infty} \mathbb{E} \bigl[G(\nu
\oplus_{r}\bar\nu_n)| \mathcal{F}_{r} \bigr] |
\mathcal{F}_{s} \Bigr]
\\
&=&\lim_{n\to\infty} \mathbb{E} \bigl[G(\nu\oplus _{r}\bar
\nu _n )|\mathcal{F}_{s} \bigr]
\\
&\ge& \mathop{\operatorname{ess}\inf}_{\bar\nu\in\mathcal{V}} \mathbb{E} \bigl[G ( \nu
\oplus_{r}\bar\nu )|\mathcal{F} _{s} \bigr]
\\
&\ge& \mathop{\operatorname{ess}\inf}_{\bar\nu\in\mathcal{V}} \mathbb{E} \bigl[G ( \nu
\oplus_{s}\bar\nu )|\mathcal{F} _{s} \bigr]
\\
&=& S^{\nu}(s),
\end{eqnarray*}
where the last inequality follows from the fact that any control $\nu
\oplus_{r}\bar\nu$, where $\bar\nu\in\mathcal{V}$, can be written in
the form $\nu\oplus_{s}(\nu\oplus_{r}\bar\nu)$; cf.
\hyperref[assVcconcat]{\textup{(C1)}}.%\vspace*{6pt}

\textit{Step} 2. There exists a family of c\`adl\`ag martingales
$\{M^{\nu}, \nu \in\mathcal{V}\} \subset\mathcal{M}_{t,p}$ such that
$S^\nu(r) \ge M^{\nu} (r)$ $\mathbb{P}$-a.s. for all $r\in[t,T]$ and
$\nu\in \mathcal{V}$.%\vspace*{6pt}

Fix $\nu\in\mathcal{V}$. By step~1, $S^\nu(\cdot)$ satisfies the
submartingale property. Therefore,
\[
S_+(r) (\omega):=\lim_{u\in(r,T]\cap\mathbb{Q},  u\to r} S^\nu (u) (\omega)
\qquad\mbox{for }0\leq r<T \quad\mbox{and}\quad S_+(T):=S^\nu(T)
\]
is well defined $\mathbb{P}$-a.s.; moreover, recalling that the
filtration $\mathbb{F}$ satisfies the usual conditions, $S_+$ is a
(right-continuous) submartingale satisfying $S_+(r)\geq S^\nu(r)$
\mbox{$\mathbb{P}$-}a.s. for all $r\in [t,T]$; cf.
\cite{DellacherieMeyer82}, Theorem~VI.2. Let $H\subset [t,T]$ be the
set of points where the function $r\mapsto\mathbb {E}[S^\nu(r)]$ is not
right-continuous. Since this function is increasing, $H$ is at most
countable. (If $H$ happens to be the empty set, then $S_+$ defines a
modification of $S^\nu$ and the Doob--Meyer decomposition of $S_+$
yields the result.) Consider the process
\[
\bar{S}(r):= S_+(r) \mathbf{1}_{H^c}(r) + S^\nu(r) \mathbf
{1}_{H}(r),\qquad r\in[t,T].
\]
The arguments (due to Lenglart) in the proof of
\cite{DellacherieMeyer82}, Theorem~10 of Appendix~1, show that
$\bar{S}$ is an \textit{optional modification} of $S^\nu$ and
$\mathbb{E}[\bar{S}(\tau)|\mathcal{F}_\sigma]\geq\bar{S}(\sigma )$ for
all $\sigma,\tau \in\mathcal{T}_{t}$ such that $\sigma\leq\tau$; that
is, $\bar {S}$ is a strong submartingale. Let
$N=N^{\mathfrak{u},\nu}_{t,z}$ be a right-continuous process of
class~(D) as in \hyperref[assIborneinf]{\textup{(I1)}}; then
$S^\nu(r)\geq N(r)$ $\mathbb {P}$-a.s. for all $r$ implies that
$S_+(r)\geq N(r)$ $\mathbb{P}$-a.s. for all $r$, and since both $S_+$
and $N$ are right-continuous, this shows that $S_+\geq N$ up to
evanescence. Recalling that $H$ is countable, we deduce that
$\bar{S}\geq N$ up to evanescence, and as $\bar{S}$ is bounded from
above by the martingale generated by $\bar{S}(T)$, we conclude that
$\bar{S}$ is of class~(D).

Now the decomposition result of Mertens~\cite{Mertens72}, Theorem~3,
yields that there exist a (true) martingale $\bar{M}$ and a
nondecreasing (not necessarily c\`adl\`ag) predictable process $\bar
{C}$ with $\bar{C}(t)=0$ such that
\[
\bar{S}=\bar{M}+\bar{C}
\]
and in view of the usual conditions, $\bar{M}$ can be chosen to be c\`
adl\`ag. We can now define $M^\nu:= \bar{M}- \bar{M}(t) + p$ on $[t,T]$
and $M^\nu(r):=p$ for $r\in[0,t)$, then $M^\nu\in\mathcal {M}_{t,p}$.
Noting that $\bar{M}(t)=\bar{S}(t)=S^\nu(t)\geq p$
by~(\ref{eqproofGDP1atteintcible}), we see that $M^\nu$ has the
required property
\[
M^\nu(r)\leq\bar{M}(r) \leq\bar{S}(r)=S^\nu(r)\qquad\mbox{$\mathbb{P}$-a.s. for all $r\in[t,T]$}.
\]

\textit{Step} 3.
Let $\tau\in\mathcal{T}_t$ have countably many values.
Then
\[
K \bigl(\tau,Z^{\mathfrak{u},\nu}_{t,z}(\tau) \bigr)\ge M^{\nu
}(
\tau) \qquad\mathbb{P}\mbox{-a.s. for all }\nu\in{\mathcal{V}}.
\]

Fix $\nu\in\mathcal{V}$ and $\varepsilon>0$, let $M^\nu$ be as in
step~2 and let $(t_i)_{i\geq1}$ be the distinct values of $\tau$. By
step~2, we have
\[
M^\nu(t_i) \le\mathop{\operatorname{ess}
\inf}_{\bar\nu\in
\mathcal{V}} \mathbb{E} \bigl[\ell \bigl(Z^{\mathfrak{u}, \nu
\oplus_{t_{i}}\bar\nu}_{t,z}(T)
\bigr)|\mathcal{F}_{t_{i}} \bigr] \qquad\mathbb{P}\mbox{-a.s.}, i
\geq1.
\]
Moreover, \hyperref[assreguIJK]{\textup{(R1)}} yields that for each
$i\ge1$, we can find a sequence $(z_{ij})_{j\ge1}\subset\mathcal{Z}$
and a Borel partition $(B_{ij})_{j\ge1}$ of $\mathcal{Z}$ such that
%
%e2.9 #&#
\begin{eqnarray}
\mathop{\operatorname{ess}\inf}_{\bar\nu\in\mathcal{V}}\mathbb {E} \bigl[\ell
\bigl(Z^{\mathfrak{u}, \nu\oplus_{t_{i}}\bar\nu}_{t,z}(T) \bigr)|\mathcal{F}_{t_{i}}
\bigr] (\omega) \le J\bigl(t_{i},z_{ij},\mathfrak{u}
[\nu\oplus_{t_i}\cdot]\bigr) (\omega
)+\varepsilon\nonumber
\\
\eqntext{\mbox{for $\mathbb{P}$-a.e. }\omega\in C_{ij}:=\bigl\{
Z^{\mathfrak{u},\nu}_{t,z}(t_{i}) \in B_{ij}\bigr\}.}
\end{eqnarray}
Since~\hyperref[assUcWithFixedNu]{\textup{(C3)}} and the definition of
$K$ in \hyperref[assesssupessinf=constant]{\textup{(Z4)}} yield that
$J(t_{i},z_{ij},\mathfrak{u}[\nu \oplus_{t_i}\cdot])\le
K(t_{i},z_{ij})$, we conclude by~\hyperref[assreguIJK]{\textup{(R1)}}
that
\[
M^\nu(t_i) (\omega) \le K(t_{i},z_{ij})
+ \varepsilon \le K\bigl(t_{i},Z^{\mathfrak{u},\nu}_{t,z}(t_{i})
(\omega) \bigr) + 2\varepsilon\qquad\mbox{for $\mathbb{P}$-a.e. }\omega\in
C_{ij}.
\]
Let $A_{i}:=\{\tau=t_i\}\in\mathcal{F}_\tau$. Then $(A_i\cap
C_{ij})_{ i,j\ge
1}$ forms a partition of $\Omega$, and the above shows that
\[
M^\nu(\tau)-2\varepsilon\le\sum_{i,j\ge1} K
\bigl(t_{i},Z^{\mathfrak
{u},\nu
}_{t,z}(t_{i}) \bigr)
\mathbf{1}_{A_{i}\cap C_{ij}} = K\bigl(\tau,Z^{\mathfrak{u},\nu}_{t,z}(\tau)
\bigr)\qquad\mathbb{P}\mbox{-a.s.}
\]
As $\varepsilon>0$ was arbitrary, the claim follows.%\vspace*{6pt}

\textit{Step} 4.
We can now prove (GDP1). Given $\tau\in
\mathcal{T}_t$, pick a sequence $(\tau_n)_{n\geq1}
\subset\mathcal{T}_t$ such that each $\tau_{n}$ has countably many
values and $\tau_n\downarrow\tau$ $\mathbb {P}$-a.s. In view of the
last statement of Lemma~\ref{lemselectionofepsoptimalstrat} below,
step~3 implies that
\[
\bigl(Z^{\mathfrak{u},\nu}_{t,z}(\tau_n),M^{\nu}(\tau
_n)-n^{-1} \bigr) \in\Lambda(\tau_n)\qquad\mathbb{P}\mbox{-a.s. for all $n\ge1$}.
\]
However, using that $Z^{\mathfrak{u},\nu}_{t,z}$ and $M^{\nu}$ are
c\`adl\`ag,
we have
\[
\bigl( \tau_n, Z^{\mathfrak{u},\nu}_{t,z}(\tau_n),
M^{\nu}(\tau_n)-n^{-1} \bigr) \to \bigl(
\tau,Z^{\mathfrak{u},\nu}_{t,z}(\tau),M^{\nu}(\tau ) \bigr)\qquad\mathbb{P}\mbox{-a.s. as }n\to\infty,
\]
so that, by the definition of $\bar\Lambda$, we deduce that
$(Z^{\mathfrak{u},\nu}_{t,z}(\tau),M^{\nu}(\tau))\in\bar\Lambda (\tau)$
$\mathbb{P}\mbox{-a.s.} $%\vspace*{6pt}

%
%le2.4 #&#
\begin{lemma}\label{lemselectionofepsoptimalstrat}
Let Assumptions \hyperref[assUcconcat]{\textup{(C2)}},
\hyperref[assAdmiA]{\textup{(C4)}},
\hyperref[assZegauxsistrategalessurevenement]{\textup{(Z1)}} and
\hyperref[assesssupessinf=constant]{\textup{(Z4)}} hold true. For each
$\varepsilon>0$, there exists a mapping $\mu^{\varepsilon }\dvtx
[0,T]\times\mathcal{Z} \to\mathfrak{U}$ such that
\[
J \bigl(t,z,\mu^{\varepsilon}(t,z) \bigr) \ge K(t,z)-\varepsilon \qquad
\mathbb{P}\mbox{-a.s. for all }(t,z)\in[0,T]\times\mathcal{Z}.
\]
In particular, if $(t,z,p)\in[0,T]\times\mathcal{Z}\times\mathbb
{R}$, then $K(t,z)> p$ implies
$(z,p)\in\Lambda(t)$.
\end{lemma}

\begin{pf}
Since $K(t,z)$ was defined
in~\hyperref[assesssupessinf=constant]{\textup{(Z4)}} as the essential
supremum of $J (t,z,\mathfrak{u})$ over $\mathfrak{u}$, there exists a
sequence $(\mathfrak{u}^{k}(t,z))_{k\ge1}\subset\mathfrak{U}$ such that
%
%e2.10 #&#
\begin{equation}
\label{eqesssup=supsurfamilledenombrable} \sup_{k\ge1} J \bigl(t,z,\mathfrak{u}^{k}(t,z)
\bigr) = K(t,z)\qquad\mathbb{P}\mbox{-a.s.}
\end{equation}
Set $\Delta_{t,z}^{0}:=\varnothing$ and define inductively the
$\mathcal{F} _t$-measurable sets
\[
\Delta_{t,z}^{k}:= \bigl\{J \bigl(t,z,\mathfrak{u}^{k}(t,z)
\bigr)\ge K(t,z)-\varepsilon \bigr\} \Bigm\backslash \bigcup
_{j=0}^{k-1} \Delta_{t,z}^{j},\qquad k\ge1.
\]
By (\ref{eqesssup=supsurfamilledenombrable}), the family $\{ \Delta
_{t,z}^{k}, k\ge1\}$ forms a partition of $\Omega$. Clearly, each
$\Delta_{t,z}^{k}$ (seen as a constant family) satisfies the
requirement of \hyperref[assAdmiA]{\textup{(C4)}} since it does not
depend on~$\nu$ and therefore belongs to $\mathfrak{F}_{t}$. Hence
after fixing some $\mathfrak {u}_0\in\mathfrak{U}$,
\hyperref[assUcconcat]{\textup{(C2)}} implies that
\[
\mu^{\varepsilon}(t,z):=\mathfrak{u}_0\oplus
_{t}\sum_{k\ge1}
\mathfrak{u}^{k}(t,z) \mathbf{1}_{\Delta
_{t,z}^{k}} \in\mathfrak{U},
\]
while \hyperref[assZegauxsistrategalessurevenement]{\textup{(Z1)}}
ensures that
\begin{eqnarray*}
J \bigl(t,z,\mu^{\varepsilon}(t,z) \bigr) &=&\mathop{\operatorname{ess}
\inf}_{\nu\in\mathcal{V}} \mathbb {E} \bigl[\ell \bigl(Z^{\mu
^{\varepsilon}(t,z),\nu
}_{t,z}(T)
\bigr)|\mathcal{F}_{t} \bigr]
\\
&=&\mathop{\operatorname{ess}\inf}_{\nu\in\mathcal{V}} \mathbb {E} \biggl[\sum
_{k\ge1}\ell \bigl(Z^{\mathfrak{u}
^{k}(t,z),\nu}_{t,z}(T)
\bigr)\mathbf{1}_{ \Delta
_{t,z}^{k}}|\mathcal {F}_{t} \biggr]
\\
&=& \mathop{\operatorname{ess}\inf}_{\nu\in\mathcal{V}} \sum
_{k\ge
1}\mathbb{E} \bigl[\ell \bigl(Z^{\mathfrak{u}^{k}(t,z),\nu
}_{t,z}(T)
\bigr)|\mathcal{F}_{t} \bigr]\mathbf{1}_{ \Delta
_{t,z}^{k}},
\end{eqnarray*}
where the last step used that $\Delta_{t,z}^{k}$ is $\mathcal
{F}_t$-measurable. Since
\[
\mathbb{E} \bigl[\ell \bigl(Z^{\mathfrak{u}^{k}(t,z),\nu
}_{t,z}(T) \bigr)|
\mathcal{F}_{t} \bigr]\ge J\bigl(t,z,\mathfrak{u}^{k}(t,z)
\bigr)
\]
by the definition of $J$, it follows by the definition of $\{\Delta
_{t,z}^{k}, k\ge1\}$ that
\[
J \bigl(t,z,\mu^{\varepsilon}(t,z) \bigr) \ge \sum
_{k\ge1} J \bigl(t,z,\mathfrak{u}^{k}(t,z) \bigr)
\mathbf{1}_{\Delta_{t,z}^{k}}\ge K(t,z)-\varepsilon\qquad\mathbb{P}\mbox{-a.s.}
\]
as required.
\end{pf}

%
%re2.5 #&#
\begin{remark}
Let us mention that the GDP could also be formulated using families of
submartingales $\{S^\nu,\nu\in\mathcal{V}\}$ rather than martingales.
Namely, in (GDP1), these would be the processes defined
by~(\ref{eqdefSnu}). However, such a formulation would not be
advantageous for applications as in Section~\ref{secexample} because we
would then need an additional control process to describe the (possibly
very irregular) finite variation part of $S^\nu$. The fact that the
martingales $\{ M^\nu,\nu\in\mathcal{V}\}$ are actually sufficient to
obtain a useful GDP can be explained heuristically as follows: the
relevant situation for the dynamic programming equation corresponds to
the adverse player choosing an (almost) optimal control $\nu$, and then
the value process $S^\nu$ will be (almost) a martingale.
\end{remark}

%s2.4 #&#
\subsection{Proof of (GDP2)}\label{secproofGDP2}

In the sequel, we fix $(t,z,p)\in[0,T]\times\mathcal{Z}\times
\mathbb{R}$ and let $\iota>0$,
$\mathfrak{u}\in\mathfrak{U}$, $\{M^{\nu},\nu\in\mathcal{V}\}\in
\mathfrak{M}
_{t,p}$, $\{\tau^{ \nu},
\nu\in\mathcal{V}\} \in\mathfrak{T}_{t}$ and $L^{\mathfrak{u},\nu}_{t,z}$
be as in (GDP2).
We shall use the dyadic discretization for the stopping times $\tau
^{\nu
}$; that is,
given $n\ge1$, we set
\[
\tau_{n}^{\nu}=\sum_{0\le i\le2^{n}-1}
t^{n}_{i+1}\mathbf{1}_{
(t^{n}_{i},t^{n}_{i+1}]}\bigl(
\tau^{\nu}\bigr)\qquad\mbox{where } t^n_i=i
2^{-n}T\mbox{ for }0\le i\le2^{n}.
\]
We shall first state the proof under the additional assumption that
%
%
%e2.11 #&#
\begin{equation}
\label{eqassumptMstopped} M^\nu(\cdot)=M^\nu\bigl(\cdot\wedge
\tau^\nu\bigr)\qquad\mbox{for all }\nu \in\mathcal{V}.
\end{equation}

\textit{Step} 1.
Fix $\varepsilon>0$ and $n\ge1$. There
exists $\mathfrak {u}^{\varepsilon}_{n} \in\mathfrak{U}$ such that
\[
\mathbb{E} \bigl[\ell \bigl(Z^{\mathfrak{u}^\varepsilon_{n},\nu
}_{t,z}(T) \bigr)|
\mathcal{F}_{\tau^{\nu}_{n}} \bigr] \ge K \bigl(\tau^{\nu}_{n},Z^{\mathfrak{u},\nu}_{t,z}
\bigl(\tau^{\nu
}_{n}\bigr) \bigr) -\varepsilon \qquad\mathbb{P}
\mbox{-a.s. for all }\nu\in\mathcal{V}.
\]

We fix $\varepsilon>0$ and $n\ge1$. It follows from
\hyperref[assreguIJK]{\textup{(R1)}} and
\hyperref[assUcconcat]{\textup{(C2)}} that, for each $i\le2^{n}$, we
can find a Borel partition $(B_{ij})_{j\ge1}$ of $\mathcal{Z}$ and a
sequence $(z_{ij})_{j\ge1}\subset\mathcal{Z}$ such that, for all $\bar
\mathfrak{u}\in\mathfrak{U}$ and $\nu\in\mathcal{V}$,
%
%e2.12 #&#
\begin{equation}\label{eqproofGDP2Step2reg1}
\mathbb{E} \bigl[\ell
\bigl(Z^{\mathfrak{u}\oplus_{t^{n}_{i}}\bar \mathfrak{u},\nu}_{t,z}(T)
\bigr)| \mathcal{F}_{t^{n}_{i}} \bigr](\omega) \ge
I\bigl(t^{n}_{i},z_{ij}, \mathfrak{u}\oplus_{t^{n}_{i}}\bar\mathfrak
{u},\nu \bigr) (\omega) - \varepsilon
\end{equation}
and
%e2.13 #&#
%e2.14 #&#
\begin{eqnarray}\label{eqproofGDP2Step2reg2}
K \bigl(t^{n}_{i},z_{ij} \bigr) \ge K
\bigl(t^{n}_{i},Z^{\mathfrak{u},\nu}_{t,z}
\bigl(t^{n}_{i}\bigr) (\omega ) \bigr) - \varepsilon
\nonumber\\[-8pt]\\[-8pt]
\eqntext{\mbox{for $\mathbb{P}$-a.e. $\omega\in C^{\nu}_{ij}:=\bigl
\{ Z^{\mathfrak{u},\nu}_{t,z}\bigl(t^{n}_{i}\bigr)\in
B_{ij}\bigr\}$.}}
\end{eqnarray}
Let $\mu^{\varepsilon}$ be as in
Lemma~\ref{lemselectionofepsoptimalstrat},
$\mathfrak{u}^{\varepsilon}_{ij}:=\mu^{\varepsilon }(t^{n}_{i},z_{ij})$
and $A^{\nu}_{ij}:=C^{\nu}_{ij}\cap\{\tau^{\nu}_{n}=t^{n}_{i}\}$, and
consider the mapping
\[
\nu\mapsto\mathfrak{u}^{\varepsilon}_{n}[\nu]:= \mathfrak{u}[\nu ]
\oplus_{\tau^{\nu}_{n}} \sum_{j\ge1, i\le2^n}
\mathfrak{u}^{\varepsilon}_{ij}[\nu]\mathbf{1}_{A^{
\nu}_{ij}}.
\]
Note that
\hyperref[assZegauxsicontroleegauxsurevenement]{\textup{(Z2)}} and
\hyperref[assAdmiA]{\textup{(C4)}} imply that
$\{C^{\nu}_{ij},\nu\in\mathcal{V}\}_{j\ge 1}\subset
\mathfrak{F}_{t^{n}_{i}}$ for each $i\le2^{n}$. Similarly, it follows
from \hyperref[assAdmitau]{\textup{(C6)}} and the definition of
$\tau^{\nu}_{n}$ that the families
$\{\{\tau^{\nu}_{n}=t^{n}_{i}\},\allowbreak \nu\in\mathcal{V}\}$ and $\{\{\tau
^{\nu}_{n}=t^{n}_{i}\}^{c},\nu\in\mathcal{V}\}$ belong to $\mathfrak{F}
_{t^{n}_{i}}$. Therefore, an induction (over $i$) based on
\hyperref[assUcconcat]{\textup{(C2)}} yields that
$\mathfrak{u}^{\varepsilon}_{n}\in\mathfrak{U}$. Using successively
(\ref{eqproofGDP2Step2reg1}),
\hyperref[assZegauxsistrategalessurevenement]{\textup{(Z1)}}, the
definition of~$J$, Lemma \ref{lemselectionofepsoptimalstrat} and
(\ref{eqproofGDP2Step2reg2}), we deduce that for $\mathbb{P}$-a.e.
$\omega\in A^{\nu}_{ij}$,
\begin{eqnarray*}
\mathbb{E} \bigl[\ell \bigl(Z^{\mathfrak{u}^\varepsilon_{n},\nu
}_{t,z}(T) \bigr)|
\mathcal{F}_{\tau^{\nu
}_{n}} \bigr](\omega) & \ge& I \bigl(t^{n}_{i},z_{ij},
\mathfrak{u}^{\varepsilon
}_{ij},\nu \bigr) (\omega) -\varepsilon
\\
& \ge& J \bigl(t^{n}_{i},z_{ij},
\mu^{\varepsilon
}\bigl(t^{n}_{i},z_{ij}\bigr)
\bigr) (\omega) - \varepsilon
\\
&\ge& K \bigl(t^{n}_{i},z_{ij} \bigr) -2
\varepsilon
\\
&\ge& K \bigl(t^{n}_{i},Z^{\mathfrak{u},\nu}_{t,z}
\bigl(t^{n}_{i}\bigr) (\omega ) \bigr) -3\varepsilon
\\
&=& K \bigl(\tau^{\nu}_{n}(\omega),Z^{\mathfrak{u},\nu}_{t,z}
\bigl(\tau ^{\nu
}_{n}\bigr) (\omega) \bigr) -3\varepsilon.
\end{eqnarray*}
As $\varepsilon>0$ was arbitrary and $\bigcup_{i,j}A^\nu_{ij}=\Omega$
$\mathbb{P}\mbox{-a.s.}$, this proves the claim. %\vspace*{6pt}

\textit{Step} 2.
Fix $\varepsilon>0$ and $n\ge1$. For all
$\nu\in\mathcal {V}$, we have
\[
\mathbb{E} \bigl[\ell \bigl(Z^{\mathfrak{u}^\varepsilon_{n},\nu
}_{t,z}(T) \bigr)|
\mathcal{F}_{\tau^{\nu
}_{n}} \bigr](\omega)\ge M^{\nu}\bigl(
\tau^{\nu}_{n}\bigr) (\omega) -\varepsilon\qquad\mbox{for $
\mathbb{P}$-a.e. $\omega\in E_{n}^{\nu}$},
\]
where
\[
E_{n}^{\nu}:= \bigl\{ \bigl(\tau_{n}^{\nu},
Z^{\mathfrak{u},\nu
}_{t,z}\bigl(\tau^{\nu
}_{n}
\bigr),M^{\nu} \bigl(\tau^{\nu}_{n}\bigr) \bigr) \in
B_{\iota} \bigl(\tau^{\nu}, Z^{\mathfrak{u},\nu}_{t,z}
\bigl(\tau^{\nu}\bigr),M^{\nu} \bigl(\tau^{\nu}\bigr)
\bigr) \bigr\}.
\]

Indeed, since $ (Z^{\mathfrak{u},\nu}_{t,z}(\tau^{\nu}),M^{\nu
}(\tau^{\nu }) )\in\Lambda_{\iota}(\tau^{\nu})$ $\mathbb
{P}\mbox{-a.s.}$, the definition of $\Lambda_{\iota}$ entails that $
(Z^{\mathfrak{u},\nu}_{t,z}(\tau^{\nu}_{n}),M^{\nu}(\tau ^{\nu}_{n})
)\in{\Lambda}(\tau^{\nu}_{n})$ for $\mathbb{P}$-a.e. $\omega\in
E_{n}^{\nu}$. This, in turn, means that
\[
K \bigl(\tau^{\nu}_{n}(\omega),Z^{\mathfrak{u},\nu}_{t,z}
\bigl(\tau ^{\nu}_{n}\bigr) (\omega ) \bigr) \ge
M^{\nu} \bigl(\tau^{\nu}_{n}\bigr) (\omega)\qquad
\mbox{for $\mathbb{P}$-a.e. $\omega\in E_{n}^{\nu}$.}
\]
Now the claim follows from step~1. [In all this, we actually have
$M^{\nu} (\tau^{\nu}_{n})=M^{\nu} (\tau^{\nu})$
by~(\ref{eqassumptMstopped}), a fact we do not use here.]%\vspace*{6pt}

\textit{Step} 3.
Let $L^{\nu}:=L^{\mathfrak{u},\nu}_{t,z}$ be the
process from \hyperref[assIbornesup]{\textup{(I2)}}. Then
\[
K(t,z) \ge p- \varepsilon- \sup_{\nu\in\mathcal{V}} \mathbb {E} \bigl[
\bigl( L^{ \nu}\bigl(\tau^{\nu
}_{n}\bigr) -
M^{\nu}\bigl(\tau^{\nu}_{n}\bigr)
\bigr)^{-}\mathbf {1}_{(E_{n}^{\nu
})^{c}} \bigr].
\]

Indeed, it follows from step 2 and
\hyperref[assIbornesup]{\textup{(I2)}} that
\begin{eqnarray*}
&&\mathbb{E} \bigl[\ell \bigl(Z^{\mathfrak{u}^\varepsilon_{n},\nu
}_{t,z}(T) \bigr)|
\mathcal{F}_{t} \bigr]
\\
&&\qquad \ge \mathbb{E} \bigl[M^{\nu}\bigl(\tau^\nu_{n}
\bigr)\mathbf {1}_{E_{n}^{\nu}}|\mathcal{F}_{t} \bigr]-\varepsilon +
\mathbb{E} \bigl[\mathbb{E} \bigl[\ell \bigl(Z^{\mathfrak
{u}^\varepsilon_{n},\nu}_{t,z}(T)
\bigr)|\mathcal{F}_{\tau^\nu
_{n}} \bigr]\mathbf{1}_{(E_{n}^{\nu})^{c}}|
\mathcal{F}_{t} \bigr]
\\
&&\qquad \ge \mathbb{E} \bigl[ M^{\nu}\bigl(\tau^{\nu}_{n}
\bigr)|\mathcal{F}_{t} \bigr] - \mathbb{E} \bigl[ M^{\nu}
\bigl(\tau^{\nu}_{n}\bigr) \mathbf{1}_{(E_{n}^{\nu})^{c}} |
\mathcal{F}_{t} \bigr] - \varepsilon+ \mathbb{E} \bigl[L^{\nu}
\bigl(\tau^{\nu}_{n}\bigr)\mathbf{1}_{(E_{n}^{\nu
})^{c}} |
\mathcal{F}_{t} \bigr]
\\
&&\qquad = p- \varepsilon+ \mathbb{E} \bigl[ \bigl( L^{\nu}\bigl(
\tau^{\nu}_{n}\bigr) - M^{\nu}\bigl(\tau
^{\nu}_{n}\bigr) \bigr)\mathbf{1}_{(E_{n}^{\nu})^{c}} |
\mathcal{F}_{t} \bigr].
\end{eqnarray*}
By the definitions of $K$ and $J$, we deduce that
\begin{eqnarray*}
K(t,z) &\ge& J\bigl(t,z,\mathfrak{u}^{\varepsilon}_{n}\bigr)
\\
&\ge& p- \varepsilon +\mathop{\operatorname{ess}\inf}_{\nu\in\mathcal{V}} \mathbb
{E} \bigl[ \bigl( L^{\nu
}\bigl(\tau^{\nu}_{n}\bigr)
- M^{\nu
}\bigl(\tau^{\nu}_{n}\bigr) \bigr)
\mathbf{1}_{(E_{n}^{\nu})^{c}} |\mathcal {F}_{t} \bigr].
\end{eqnarray*}
Since $K$ is deterministic, we can take expectations on both sides to
obtain that
%
%e2.15 #&#
\begin{eqnarray}
K(t,z) \ge p- \varepsilon +\mathbb{E} \bigl[\mathop{\operatorname{ess}
\inf}_{\nu\in\mathcal
{V}} \mathbb{E} \bigl[ Y^{\nu} |\mathcal{F}_{t}
\bigr] \bigr],\nonumber
\\
\eqntext{\mbox{where }Y^{\nu}:= \bigl(L^{\nu}\bigl(
\tau^{\nu}_{n}\bigr) - M^{\nu}\bigl(
\tau^{\nu}_{n}\bigr) \bigr)\mathbf{1}_{(E_{n}^{\nu})^{c}}.}
\end{eqnarray}
The family $\{ \mathbb{E} [ Y^{\nu} |\mathcal{F}_{t}
],\nu\in\mathcal{V}\}$ is directed downward; to see this, use
\hyperref[assVcconcat]{\textup{(C1)}},
\hyperref[assZegauxsicontroleegauxsurevenement]{\textup{(Z2)}},
\hyperref[assMegauxsicontroleegauxsurevenement]{\textup{(Z3)}},
\hyperref[asstauegauxsicontroleegauxsurevenement]{\textup{(C5)}} and
the last statement in \hyperref[assIbornesup]{\textup{(I2)}}, and argue
as in step~1 of the proof of~(GDP1) in Section~\ref{secproofGDP1}. It
then follows that we can find a sequence $(\nu_{k})_{k\ge1}\subset
\mathcal{V}$ such that $\mathbb{E} [ Y^{\nu_{k}} |\mathcal{F}_{t} ]$
decreases $\mathbb{P}\mbox{-a.s.}$ to $\mathop{\operatorname {ess}\inf}
_{\nu\in\mathcal{V}} \mathbb{E} [ Y^{\nu} |\mathcal {F}_{t} ]$ (cf.
\cite{Neveu75}, Proposition~VI-1-1) so that the claim follows by
monotone convergence.%\vspace*{6pt}

\textit{Step} 4. We have
\[
\lim_{n\to\infty}\sup_{\nu\in\mathcal{V}} \mathbb{E} \bigl[
\bigl( L^{ \nu}\bigl(\tau^{\nu
}_{n}\bigr) -
M^{\nu}\bigl(\tau^{\nu}_{n}\bigr)
\bigr)^{-}1_{(E_{n}^{\nu})^{c}} \bigr]= 0 \qquad\mathbb{P}\mbox{-a.s.}
\]

Indeed, since $M^{\nu} (\tau^{\nu}_{n})=M^{\nu} (\tau^{\nu})$
by~(\ref{eqassumptMstopped}), the uniform integrability assumptions in
Theorem~\ref{thmGDP} yield that $\{ ( L^{ \nu}(\tau^{\nu }_{n}) -
M^{\nu}(\tau^{\nu}_{n})  )^{-}\dvtx  n\geq1,\nu\in\mathcal {V}\}$ is
again uniformly integrable. Therefore, it suffices to prove that\break
$\sup_{\nu \in\mathcal{V}}\mathbb{P} \{{(E_{n}^{\nu})^{c}} \} \to 0$.
To see this, note that for $n$ large enough, we have
$|\tau^{\nu}_{n}-\tau^{\nu}|\le2^{-n}T\le \iota /2$ and hence
\[
\mathbb{P} \bigl\{{\bigl(E_{n}^{\nu}\bigr)^{c}}
\bigr\} \le\mathbb{P} \bigl\{ {d_{\mathcal{Z}} \bigl(Z^{\mathfrak{u},\nu}_{t,z}
\bigl(\tau^{\nu
}_{n}\bigr),Z^{\mathfrak{u},\nu}_{t,z}
\bigl(\tau^{\nu}\bigr) \bigr)\ge\iota/2 } \bigr\},
\]
where we have used that $M^{\nu} (\tau^{\nu}_{n})=M^{\nu} (\tau
^{\nu})$. Using once more that $|\tau^{\nu}_{n}-\tau^{\nu}|\le2^{-n}T$,
the claim then follows from
\hyperref[assreguprobaZ]{\textup{(R2)}}.%\vspace*{6pt}

\textit{Step} 5. The additional assumption~(\ref{eqassumptMstopped})
entails no loss of generality.%\vspace*{6pt}

Indeed, let $\tilde{M}^\nu$ be the stopped martingale $M^\nu(\cdot
\wedge \tau^\nu)$. Then $\{\tilde{M}^\nu,\nu\in\mathcal{V}\}\subset
\mathcal{M}_{t,p}$. Moreover, since
$\{M^\nu,\nu\in\mathcal{V}\}\in\mathfrak {M}_{t,p}$ and $\{
\tau^\nu,\nu\in \mathcal{V}\}\in\mathfrak{T}_t$, we see
from~\hyperref[assMegauxsicontroleegauxsurevenement]{\textup{(Z3)}}
and~\hyperref[asstauegauxsicontroleegauxsurevenement]{\textup{(C5)}}
that $\{ \tilde {M}^\nu,\nu\in\mathcal{V}\}$ again satisfies the
property stated in
\hyperref[assMegauxsicontroleegauxsurevenement]{\textup{(Z3)}}.
Finally, we have that the set $ \{\tilde{M}^{\nu}(\tau^\nu)^{+}\dvtx
\nu\in\mathcal{V} \} $ is uniformly integrable like $
\{M^{\nu}(\tau^\nu)^{+}\dvtx \nu\in\mathcal {V} \}$, since these sets
coincide. Hence, $\{\tilde{M}^\nu,\nu\in\mathcal{V}\}$ satisfies all
properties required in (GDP2), and of course
also~(\ref{eqassumptMstopped}). To be precise, it is not necessarily
the case that
$\{\tilde{M}^\nu,\nu\in\mathcal{V}\}\in\mathfrak{M}_{t,p}$; in fact, we
have made no assumption whatsoever about the richness of
$\mathfrak{M}_{t,p}$. However, the previous properties are all we have
used in this proof and hence, we may indeed replace $M^\nu$ by
$\tilde{M}^\nu$ for the purpose of proving (GDP2).

We can now complete the proof of (GDP2): in view of step~4, step~3
yields that $K(t,z)\geq p-\varepsilon$, which by
Lemma~\ref{lemselectionofepsoptimalstrat} implies the assertion that
$(z,p-\varepsilon)\in \Lambda (t)$.

%s2.5 #&#
\subsection{Proof of Corollary~\protect\texorpdfstring{\ref{corGDPwithlessregularity}}{2.3}}\label{secproofcorGDPwithlessregularity}

\textit{Step} 1. Assume that $\ell$ is bounded and Lipschitz
continuous. Then \ref{assI} and \hyperref[assreguIJK]{\textup{(R1)}}
are satisfied.%\vspace*{6pt}

Assumption \ref{assI} is trivially satisfied; we prove that
(\ref{eqZlipschitz}) implies Assumption~\hyperref[assreguIJK]{\textup{(R1)}}. Let $t\leq s\leq T$ and
$(\mathfrak{u},\nu)\in\mathfrak{U}\times \mathcal {V}$. Let $c$ be the
Lipschitz constant of~$\ell$. By~(\ref{eqZlipschitz}), we have
%
%
%e2.16 #&#
\begin{eqnarray}\label{eqlipEstimate}
\bigl\llvert \mathbb{E} \bigl[\ell \bigl( Z^{\mathfrak{u},\nu}_{t,z}(T)
\bigr)-\ell \bigl( Z^{\mathfrak{u},\nu}_{s,z'}(T) \bigr)|
\mathcal{F}_s \bigr] \bigr\rrvert &\le& c\mathbb{E} \bigl[\bigl\llvert
Z^{\mathfrak{u},\nu
}_{t,z}(T)-Z^{\mathfrak{u},\nu}_{s,z'}(T)\bigr
\rrvert |\mathcal{F} _s \bigr]
\nonumber
\\
&\le& cC \bigl\llvert Z^{\mathfrak{u},\nu}_{t,z}(s)-z'\bigr
\rrvert
\end{eqnarray}
for all $z,z'\in\mathbb{R}^d$. Let $(B_j)_{j\ge1}$ be any Borel
partition of $\mathbb{R}
^d$ such that the diameter of $B_j$ is less than $\varepsilon/(cC)$,
and let $z_j\in B_j$ for each $j\ge1$. Then
\[
\bigl\llvert \mathbb{E} \bigl[\ell \bigl( Z^{\mathfrak{u},\nu}_{t,z}(T)
\bigr)-\ell \bigl( Z^{\mathfrak{u},\nu}_{s,z_j}(T) \bigr)|
\mathcal{F}_s \bigr] \bigr\rrvert \le\varepsilon\qquad\mbox{on }
C^{\mathfrak{u},\nu}_j:= \bigl\{ Z^{\mathfrak{u},\nu
}_{t,z}(s)
\in B_j \bigr\},
\]
which implies the first property in
\hyperref[assreguIJK]{\textup{(R1)}}. In particular, let
$\bar\nu\in\mathcal{V}$, then using
\hyperref[assVcconcat]{\textup{(C1)}}, we have
\[
\bigl\llvert \mathbb{E} \bigl[\ell \bigl( Z^{\mathfrak{u},\nu\oplus_s\bar\nu
}_{t,z}(T)
\bigr)-\ell \bigl( Z^{\mathfrak{u},\nu\oplus_s\bar\nu}_{s,z_j}(T) \bigr)|
\mathcal{F}_s \bigr] \bigr\rrvert \le\varepsilon\qquad\mbox{on }
C^{\mathfrak{u},\nu\oplus_s\bar
\nu}_j.
\]
Since $C_j^{\mathfrak{u},\nu\oplus_s\bar\nu}=C^{\mathfrak{u},\nu }_j$
by \hyperref[assZegauxsicontroleegauxsurevenement]{\textup{(Z2)}}, we
may take the essential infimum over $\bar\nu\in\mathcal{V}$ to conclude
that
\[
\mathop{\operatorname{ess}\inf}_{\bar\nu\in\mathcal{V}}\mathbb {E} \bigl[\ell \bigl(
Z^{\mathfrak{u},\nu\oplus_s\bar\nu
}_{t,z}(T) \bigr)|\mathcal{F}_s \bigr] \le J
\bigl(s,z_j,\mathfrak{u} [\nu\oplus _s\cdot ] \bigr)+
\varepsilon \qquad\mbox{on } C^{\mathfrak{u},\nu}_j,
\]
which is the second property in \hyperref[assreguIJK]{\textup{(R1)}}.
Finally, the last property in~\hyperref[assreguIJK]{\textup{(R1)}} is a
direct consequence of~(\ref{eqlipEstimate}) applied with
$t=s$.%\vspace*{6pt}

\textit{Step} 2. We now prove the corollary under the additional
assumption that \mbox{$|\ell(z)|\leq C$}; we shall reduce to the Lipschitz
case by inf-convolution. Indeed, if we define the functions $\ell_k$ by
\[
\ell_k(z)=\inf_{z'\in\mathbb{R}^d}\bigl\{ \ell
\bigl(z'\bigr) + k\bigl|z'-z\bigr|\bigr\},\qquad k\geq1
\]
then $\ell_k$ is Lipschitz continuous with Lipschitz constant $k$,
$|\ell_k|\leq C$, and $(\ell_k)_{k\ge1}$ converges pointwise to $\ell$.
Since $\ell$ is continuous and the sequence $(\ell_k)_{k\ge1}$ is
monotone increasing, the convergence is uniform on compact sets by
Dini's lemma. That is, for all $n\geq1$,
%
%e2.17 #&#
\begin{equation}
\label{eqellketellprochessurcompact} \sup_{z \in\mathbb{R}^d,  |z |\le n} \bigl\llvert \ell_k(z)
- \ell (z)\bigr\rrvert \le \epsilon^n_k,
\end{equation}
where $(\epsilon^n_k)_{k\ge1}$ is a sequence of numbers such that
$\lim_{k\to\infty}\epsilon^n_k=0$. Moreover, (\ref{eqassZgrowth})
combined with Chebyshev's inequality imply that
%
%e2.18 #&#
\begin{equation}
\label{eqprobaZplusgrandquentendvers0}
\mathop{\operatorname{ess}\sup}_{(\mathfrak{u},\nu)\in\mathfrak
{U}\times\mathcal{V}} \mathbb {P} \bigl
\{{\bigl|Z^{\mathfrak{u},\nu}_{t,z}(T)\bigr|\ge n|\mathcal {F}_{t}} \bigr
\} \le\bigl(\varrho(z)/n\bigr)^{\bar{q}}.
\end{equation}
Combining (\ref{eqellketellprochessurcompact}) and
(\ref{eqprobaZplusgrandquentendvers0}) and using the fact that $\ell
_{k}- \ell$ is bounded by $2C$ then leads to
%
%e2.19 #&#
\begin{equation}
\label{eqdifferencelkatltendvers0} \quad\mathop{\operatorname{ess}\sup}_{(\mathfrak{u},\nu)\in\mathfrak
{U}\times\mathcal{V}} \mathbb{E} \bigl[
\bigl\llvert \ell_k \bigl(Z^{\mathfrak{u},\nu
}_{t,z}(T)
\bigr)-\ell \bigl(Z^{\mathfrak{u},\nu
}_{t,z}(T) \bigr)\bigr\rrvert |
\mathcal{F}_{t} \bigr] \le\epsilon^n_k + 2C
\bigl(\varrho(z)/n\bigr)^{\bar{q}}.
\end{equation}
Let $O$ be a bounded subset of $\mathbb{R}^d$, let $\eta>0$ and let
%
%e2.20 #&#
\begin{equation}
\label{eqdefIk} I_k(t,z,\mathfrak{u},\nu)=\mathbb{E} \bigl[
\ell_k \bigl( Z^{\mathfrak{u},\nu}_{t,z}(T) \bigr)|
\mathcal{F}_{t} \bigr].
\end{equation}
Then\vspace*{-2pt} we can choose an integer $n^{\eta}_{O}$ such that $2C
(\varrho(z)/n^{\eta}_{O})^{\bar{q}}\leq\eta/2$ for all\allowbreak $z\in O$ and another integer $k^{\eta}_{O}$ such that $\epsilon^{n^{\eta
}_{O}}_{k^{\eta}_{O}}\le\eta/2$. Under\vspace*{-3pt} these conditions,
(\ref{eqdifferencelkatltendvers0})~applied to $n=n^{\eta}_{O}$ yields
that
%
%e2.21 #&#
\begin{equation}
\label{eqIkeps} \qquad\quad\mathop{\operatorname{ess}\sup}_{(\mathfrak{u},\nu)\in\mathfrak
{U}\times\mathcal{V}}\bigl\llvert
I_{k_O^{\eta}} (t,z,\mathfrak{u},\nu) - I(t,z,\mathfrak{u},\nu) \bigr\rrvert
\le\eta\qquad \mbox{for } (t,z) \in [0,T]\times O.
\end{equation}
In the sequel, we fix $(t,z,p)\in[0,T]\times\mathbb{R}^d \times
\mathbb{R}$ and a bounded set $O\subset\mathbb{R}^d$ containing $z$,
and define $J_{k^{\eta}_{O}}$, $\Lambda_{k^{\eta}_{O}}$,
$\Lambda_{k^{\eta}_{O},\iota }$ and $\bar\Lambda_{k^{\eta}_{O}}$ in
terms of $\ell_{k^{\eta}_{O}}$ instead of~$\ell$.

We now prove (GDP1$'$). To this end, suppose that $(z,p+ 2\eta)\in
\Lambda(t)$. Then (\ref{eqIkeps}) implies that $(z,p+\eta)\in
\Lambda_{k^{\eta}_{O}}(t)$. In view of step~1, we may apply (GDP1) with
the loss function $\ell_{k^{\eta}_{O}}$ to obtain $\mathfrak{u}\in
\mathfrak{U}$ and $\{
M^{\nu},\nu\in{\mathcal{V}}\}\subset\mathcal{M}_{t,p}$ such that
\[
\bigl(Z^{\mathfrak{u},\nu}_{t,z}(\tau),M^{\nu}(\tau)+\eta \bigr)
\in\bar\Lambda_{k^{\eta}_{O}}(\tau) \qquad\mathbb{P}\mbox {-a.s. for all $\nu \in\mathcal{V}$ and $\tau\in\mathcal{T}_{t}$.}
\]
Using once more (\ref{eqIkeps}), we deduce that
%
%e2.22 #&#
\begin{eqnarray}
\hspace*{-30pt}\bigl(Z^{\mathfrak{u},\nu}_{t,z}(\tau),M^{\nu}(\tau) \bigr) \in
\bar\Lambda(\tau)\nonumber
\\
\eqntext{\mathbb{P}\mbox{-a.s. for all $\nu\in\mathcal{V}$
and $\tau\in\mathcal{T}_{t}$ such that $Z^{\mathfrak{u},\nu
}_{t,z}(\tau)\in O$.}}
\end{eqnarray}
Recalling that $\{Z^{\mathfrak{u},\nu}_{t,z} (\tau^{\mathfrak {u},\nu}
),(\mathfrak{u},\nu )\in\mathfrak{U}\times\mathcal{V}\}$ is uniformly
bounded and enlarging~$O$ if necessary, we deduce that (GDP1$'$) holds
for $\ell$. [The last two arguments are superfluous as
$\ell\ge\ell_{k^{\eta}_{O}}$ already implies
$\bar\Lambda_{k^{\eta}_{O}}(\tau)\subset\bar\Lambda(\tau )$; however,
we would like to refer to this proof in a similar situation below where
there is no monotonicity.]

It remains to prove (GDP2$'$). To this end, let $\iota>0$, $\mathfrak
{u}\in\mathfrak{U} $,
$\{M^{\nu},\nu\in\mathcal{V}\}\in\mathfrak{M}_{t,p}$ and $\{ \tau^{ \nu
}, \nu\in\mathcal{V} \} \in\mathfrak{T}_{t}$ be such that
\[
\bigl(Z^{\mathfrak{u},\nu}_{t,z}\bigl(\tau^{\nu}
\bigr),M^{\nu}\bigl(\tau^{\nu
}\bigr) \bigr) \in\Lambda_{2\iota}\bigl(\tau^{\nu}\bigr) \qquad\mathbb{P}\mbox{-a.s. for all } \nu\in\mathcal{V}.
\]
For $\eta<\iota/2$, we then have
%
%e2.23 #&#
\begin{equation}\label{eqreplacebarM}
\bigl(Z^{\mathfrak{u},\nu}_{t,z}\bigl(\tau^{\nu}
\bigr),M^{\nu}\bigl(\tau^{\nu
}\bigr)+2\eta \bigr) \in\Lambda_{\iota}\bigl(\tau^{\nu}\bigr) \qquad\mathbb{P}\mbox{-a.s. for all } \nu\in\mathcal{V}.
\end{equation}
Let $\tilde{M}^\nu:=M^{\nu}+\eta$. Since $\{Z^{\mathfrak{u},\nu }_{t,z}
(\tau ^{\nu}  ),\nu\in\mathcal{V}\}$ is uniformly bounded in
$L^{\infty}$, we may assume, by enlarging\vadjust{\goodbreak} $O$ if necessary, that
$B_{\iota }(Z^{\mathfrak{u},\nu }_{t,z}(\tau^{\nu})) \subset O$
$\mathbb{P}\mbox{-a.s.}$ for all $\nu\in\mathcal{V}$. Then
(\ref{eqIkeps}) and (\ref{eqreplacebarM}) imply that
\[
\bigl(Z^{\mathfrak{u},\nu}_{t,z}\bigl(\tau^{\nu}\bigr),
\tilde{M}^{\nu
}\bigl(\tau^{\nu}\bigr) \bigr) \in\Lambda_{k^{\eta}_{O},\iota}\bigl(\tau^{\nu}\bigr) \qquad\mathbb{P}
\mbox{-a.s. for all } \nu\in\mathcal{V}.
\]
Moreover, as $\ell\le C$, (\ref{eqreplacebarM}) implies that $\tilde
{M}^\nu(\tau^{\nu})\le C$; in particular, $\{\tilde{M}^\nu(\tau ^\nu
)^+,\allowbreak \nu\in\mathcal{V}\}$ is uniformly integrable. Furthermore, as
$\ell\ge-C$, we can take $L^{\mathfrak{u},\nu}_{t,z}:=-C$
for~\hyperref[assIbornesup]{\textup{(I2)}}. In view of step~1, (GDP2)
applied with the loss function $\ell _{k^{\eta }_{O}}$ then yields that
%
%
%e2.24 #&#
\begin{equation}
\label{eqproofWithTildeM} (z,p+\eta-\varepsilon)\in\Lambda_{k^{\eta}_{O}}(t)\qquad\mbox{for
all }\varepsilon>0.
\end{equation}
To be precise, this conclusion would require that $\{\tilde{M}^\nu,\nu
\in\mathcal{V}\}\in\mathfrak{M}_{t,p+\eta}$, which is not necessarily
the case under our assumptions. However, since
$\{M^{\nu},\nu\in\mathcal{V}\}\in \mathfrak{M} _{t,p}$, it is clear
that $\{\tilde{M}^\nu,\nu\in\mathcal{V}\}$ satisfies the property
stated in
\hyperref[assMegauxsicontroleegauxsurevenement]{\textup{(Z3)}}, so
that, as in step~5 of the proof of (GDP2), there is no loss of
generality in assuming that
$\{\tilde{M}^\nu,\nu\in\mathcal{V}\}\in\mathfrak{M} _{t,p+\eta}$. We
conclude by noting that~(\ref{eqIkeps}) and (\ref{eqproofWithTildeM})
imply that $(z,p-\varepsilon)\in\Lambda(t)$ for all
$\varepsilon>0$.%\vspace*{6pt}

\textit{Step} 3. We turn to the general case. For $k\ge1$, we now
define $\ell _{k}:= (\ell\wedge k)\vee(-k)$, while $I_k$ is again
defined as in~(\ref{eqdefIk}). We also set
\[
n_k=\max \bigl\{m\geq0\dvtx  B_m(0)\subset\{\ell=
\ell_k\} \bigr\}\wedge k
\]
and note that the continuity of $\ell$ guarantees that $\lim_{k\to
\infty
} n_k=\infty$.
Given a bounded set $O\subset\mathbb{R}^d$ and $\eta>0$, we claim that
%
%e2.25 #&#
%e2.26 #&#
\begin{eqnarray}
 \mathop{\operatorname{ess}\sup}_{(\mathfrak{u},\nu)\in\mathfrak
{U}\times\mathcal{V}} \bigl\llvert
I_{k_O^{\eta}} (t,z,\mathfrak{u},\nu) - I(t,z,\mathfrak{u},\nu) \bigr\rrvert
\le\eta
\nonumber\\[-10pt]\label{eqclaimProofCorStep3} \\[-10pt]
\eqntext{\mbox{for all } (t,z) \in [0,T]\times O}
\end{eqnarray}
for any large enough integer $k_O^{\eta}$.
Indeed, let $(\mathfrak{u},\nu)\in\mathfrak{U}\times\mathcal{V}$; then
\begin{eqnarray*}
&& \bigl\llvert I_{k} (t,z,\mathfrak{u},\nu) - I(t,z,\mathfrak{u},\nu)
\bigr\rrvert
\\
&&\qquad \le \mathbb{E} \bigl[|\ell-\ell_k|
\bigl(Z^{\mathfrak{u},\nu
}_{t,z}(T) \bigr)|\mathcal{F}_{t}
\bigr]
\\
&&\qquad =  \mathbb{E} \bigl[|\ell-\ell_k| \bigl(Z^{\mathfrak{u},\nu
}_{t,z}(T)
\bigr)\mathbf{1}_{ \{Z^{\mathfrak{u},\nu
}_{t,z}(T)\notin\{\ell=\ell_k\} \}}|\mathcal{F}_{t} \bigr]
\\
&&\qquad \leq \mathbb{E} \bigl[\bigl\llvert \ell \bigl(Z^{\mathfrak{u},\nu
}_{t,z}(T)
\bigr)\bigr\rrvert \mathbf{1}_{ \{|Z^{\mathfrak{u},\nu
}_{t,z}(T)|>n_k \}}|\mathcal{F}_{t} \bigr]
\\
&&\qquad \leq C\mathbb{E} \bigl[ \bigl(1+\bigl\llvert Z^{\mathfrak{u},\nu
}_{t,z}(T)
\bigr\rrvert ^q \bigr) \mathbf{1}_{ \{|Z^{\mathfrak{u},\nu
}_{t,z}(T)|>n_k \}}|
\mathcal{F}_{t} \bigr]
\end{eqnarray*}
by~(\ref{eqassellZpoly}). We may assume that $q>0$, as otherwise we are
in the setting of step~2. Pick $\delta>0$ such that $q(1+\delta
)=\bar{q}$. Then H\"older's inequality and (\ref{eqassZgrowth}) yield
that
\begin{eqnarray*}
&&\mathbb{E} \bigl[\bigl\llvert \bigl(Z^{\mathfrak{u},\nu}_{t,z}(T) \bigr)
\bigr\rrvert ^q\mathbf{1}_{ \{|Z^{\mathfrak{u},\nu
}_{t,z}(T)|>n_k \}}|\mathcal{F}_{t}
\bigr]
\\
&&\qquad \leq \mathbb{E} \bigl[\bigl\llvert \bigl(Z^{\mathfrak{u},\nu
}_{t,z}(T)\bigr)\bigr\rrvert ^{\bar
{q}}|\mathcal{F}_{t}
\bigr]^{1/(1+\delta)} \mathbb{P} \bigl\{ \bigl|Z^{\mathfrak{u},\nu}_{t,z}(T)\bigr|>n_k|
\mathcal{F} _{t} \bigr\} ^{\delta/(1+\delta)}
\\
&&\qquad \leq \rho(z)^{\bar{q}/(1+\delta)}\bigl(\rho(z)/n_k
\bigr)^{\bar{q}\delta/(1+\delta)}.
\end{eqnarray*}
Since $\rho$ is locally bounded and $\lim_{k\to\infty} n_k=\infty$,
claim~(\ref{eqclaimProofCorStep3}) follows. We can then obtain
(GDP1$'$) and (GDP2$'$) by reducing to the result of step~2, using the
same arguments as in the proof of step~2.

%s3 #&#
\section{The PDE in the case of a controlled SDE}\label{secexample}

In this section, we illustrate how our GDP can be used to derive a
dynamic programming equation and how its assumptions can be verified in
a typical setup. To this end, we focus on the case where the state
process is determined by a stochastic differential equation with
controlled coefficients; however, other examples could be treated similarly.

%s3.1 #&#
\subsection{Setup}\label{secBMcasemainassumptions}

Let $\Omega=C([0,T];\mathbb{R}^d)$ be the canonical space of continuous
paths equipped with the Wiener measure $\mathbb{P}$, let $\mathbb
{F}=(\mathcal{F}_{t})_{t\le T}$ be the $\mathbb{P}$-augmentation of the
filtration generated by the coordinate-mapping process $W$ and let
$\mathcal{F}=\mathcal{F}_T$. We define $\mathcal{V}$, the set of
adverse controls, to be the set of all progressively measurable
processes with values in a compact subset $V$ of $\mathbb{R}^d$.
Similarly, $\mathcal{U}$ is the set of all progressively measurable
processes with values in a compact $U\subset\mathbb{R}^d$. Finally, the
set of strategies $\mathfrak{U}$ consists of all mappings
$\mathfrak{u}\dvtx \mathcal {V}\rightarrow \mathcal{U}$ which are
nonanticipating in the sense that
\[
\{\nu_1=_{(0,s]}\nu_2\} \subset \bigl\{
\mathfrak{u}[\nu _1]=_{(0,s]}\mathfrak{u}[\nu_2]
\bigr\} \qquad\mbox{for all $\nu_1,\nu_2 \in\mathcal{V}$
and $s\leq T$.}
\]
Given $(t,z)\in[0,T]\times\mathbb{R}^d$ and $(\mathfrak{u},\nu
)\in\mathfrak{U}\times\mathcal{V}$, we let
$Z^{\mathfrak{u},\nu}_{t,z}$ be the unique strong solution of the
controlled SDE
%
%e3.1 #&#
%e3.2 #&#
\begin{eqnarray}
\qquad Z(s)=z+\int_t^s \mu\bigl(
Z(r), \mathfrak{u}[\nu]_r,\nu_r\bigr) \,dr + \int
_t^s \sigma\bigl( Z(r), \mathfrak{u}[
\nu]_r,\nu_r\bigr) \,dW_r,
\nonumber\\[-12pt]\label{eqedsz} \\[-8pt]
\eqntext{s \in[t,T],}
\end{eqnarray}
where the coefficients
\[
\mu\dvtx \mathbb{R}^{d}\times U \times V\to\mathbb{R}^{d},
\qquad\sigma\dvtx \mathbb{R}^{d}\times U\times V\to\mathbb{R}
^{d\times d}
\]
are assumed to be jointly continuous in all three variables, Lipschitz
continuous with linear growth in the first variable, uniformly in the
last two and Lipschitz continuous in the second variable, locally
uniformly in the two others. Throughout this section, we assume that
$\ell\dvtx \mathbb{R}^d\to\mathbb{R}$ is a continuous function of
polynomial growth; that is, (\ref{eqassellZpoly}) holds true for some
constants $C$ and $q$. Since $Z^{\mathfrak{u},\nu}_{t,z}(T)$ has~moments of all orders, this implies that the finiteness
condition~(\ref{eqelldansL1}) is satisfied.

In view of the martingale representation theorem, we can identify the
set $\mathcal{M}_{t,p}$ of martingales with the set $\mathcal{A}$ of
all progressively measurable $d$-dimensional processes $\alpha$ such
that $\int\alpha \,dW$ is a (true) martingale. Indeed, we have
$\mathcal{M}_{t,p}=\{ P^\alpha _{t,p},\alpha\in\mathcal{A}\}$, where
\[
P_{t,p}^\alpha(\cdot)=p+\int_t^\cdot
\alpha_s \,dW_s.
\]
We shall denote by $\mathfrak{A}$ the set of all mappings
$\mathfrak{a}[\cdot]$: $\mathcal{V} \mapsto\mathcal{A}$ such that
\[
\{\nu_1=_{(0,s]}\nu_2\} \subset \bigl\{
\mathfrak{a}[\nu _1]=_{(0,s]}\mathfrak{a}[\nu_2]
\bigr\} \qquad\mbox{for all $\nu_1,\nu_2 \in\mathcal{V}$
and $s\leq T$}.
\]
The set of all families $\{P^{\mathfrak{a}[\nu]}_{t,p},\nu\in
\mathcal{V}\}$ with $\mathfrak{a}\in \mathfrak{A}$ then forms the set
$\mathfrak{M}_{t,p}$, for any given $(t,p)\in [0,T]\times\mathbb{R}$.
Furthermore, $\mathfrak{T}_{t}$ consists of all families
$\{\tau^\nu,\nu\in \mathcal{V}\} \subset\mathcal{T}_{t}$ such that, for
some $(z,p)\in\mathbb {R}^d\times\mathbb{R}$,
$(\mathfrak{u},\mathfrak{a})\in \mathfrak{U}\times\mathfrak{A}$ and
some Borel set $O\subset [0,T]\times\mathbb {R}^d \times\mathbb{R}$,
\[
\mbox{$\tau^{\nu}$ is the first exit time of $ \bigl(\cdot,
Z^{\mathfrak{u},\nu
}_{t,z}, P^{\mathfrak{a}[\nu]}_{t,p} \bigr)$ from
$O$\qquad for all $\nu \in\mathcal{V}$.}
\]
(This includes the deterministic times $s\in[t,T]$ by the choice
$O=[0,s]\times\mathbb{R}^d\times\mathbb{R}$.) Finally, $\mathfrak{F}_{t}$
consists of all
families $\{A^\nu,\nu\in\mathcal{V}\}\subset\mathcal{F}_t$ such that
\[
A^{\nu_1}\cap\{\nu_1=_{{(0,t]}} \nu_2\} =
A^{\nu_2}\cap\{\nu_1=_{{
(0,t]}} \nu_2\}\qquad\mbox{for all }\nu_1,\nu_2\in\mathcal{V}.
\]

%
%pr3.1 #&#
\begin{proposition}\label{propverifassumptionmarkovsde}
The conditions of Corollary \ref{corGDPwithlessregularity} are
satisfied in the present setup.
\end{proposition}

\begin{pf}
The above definitions readily yield that Assumptions \ref{assC}
and\allowbreak
\hyperref[assZegauxsistrategalessurevenement]{\textup{(Z1)}}--\hyperref[assMegauxsicontroleegauxsurevenement]{\textup{(Z3)}}
are satisfied. Moreover,
Assumption~\hyperref[assesssupessinf=constant]{\textup{(Z4)}} can be
verified exactly as in~\cite{BuckdahnLi08}, Proposition~3.3. Fix any
$\bar{q}>q\vee2$; then~(\ref{eqassZgrowth}) can be obtained as follows.
Let $(\mathfrak{u},\nu ) \in\mathfrak{U} \times\mathcal{V}$ and
$A\in\mathcal{F}_t$ be arbitrary. Using the Burkholder--Davis--Gundy
inequalities, the boundedness of $U$ and $V$ and the assumptions on
$\mu$ and $\sigma$, we obtain that
\[
\mathbb{E} \Bigl[\sup_{t\le s\le\tau}\bigl\llvert Z^{\mathfrak{u},\nu
}_{t,z}(s)
\bigr\rrvert ^{\bar
{q}}\mathbf{1}_A \Bigr] \le c\mathbb{E}
\biggl[ \mathbf{1}_A+\llvert z \rrvert ^{\bar{q}}
\mathbf{1}_A + \int_t^{\tau} \sup
_{t\le s\le r}\bigl\llvert Z^{\mathfrak{u},\nu
}_{t,z}(s)\bigr
\rrvert ^{\bar{q}}\mathbf{1}_A \,dr \biggr],
\]
where $c$ is a universal constant, and $\tau$ is any stopping time such
that \mbox{$Z^{\mathfrak{u},\nu}_{t,z}(\cdot\wedge\tau)$} is bounded.
Applying Gronwall's inequality and letting $\tau\to T$, we deduce that
\[
\mathbb{E} \bigl[\bigl\llvert Z^{\mathfrak{u},\nu}_{t,z}(T) \bigr\rrvert
^{\bar{q}}{\bf1}_A \bigr] \le\mathbb{E} \Bigl[\sup
_{t\le u\le T}\bigl\llvert Z^{\mathfrak{u},\nu
}_{t,z}(u) \bigr
\rrvert ^{\bar
{q}}{\bf1}_A \Bigr] \le c\mathbb{E} \bigl[
\bigl(1+\llvert z \rrvert ^{\bar{q}}\bigr){\bf1}_A \bigr].
\]
Since $A\in\mathcal{F}_t$ was arbitrary, this
implies~(\ref{eqassZgrowth}). To verify condition~(\ref{eqZlipschitz}),
we note that the flow property yields
\[
\mathbb{E} \bigl[\bigl\llvert Z^{\mathfrak{u}\oplus_s\bar\mathfrak
{u},\nu\oplus_s\bar\nu}_{t,z}(T)-Z^{\bar
\mathfrak{u},\nu\oplus_s\bar\nu}_{s,z'}(T)
\bigr\rrvert {\bf1}_A \bigr] = \mathbb{E} \bigl[\bigl\llvert
Z^{\bar\mathfrak{u},\nu\oplus_s\bar\nu
}_{s,Z^{\mathfrak{u},\nu
}_{t,z}(s)}(T)-Z^{\bar\mathfrak{u},\nu\oplus_s\bar\nu}_{s,z'}(T) \bigr
\rrvert {\bf1}_A \bigr]
\]
and estimate the right-hand side with the above arguments. Finally, the
same arguments can be used to
verify~\hyperref[assreguprobaZ]{\textup{(R2)}}.
\end{pf}

%

%
%re3.2 #&#
\begin{remark}
We emphasize that our definition of a strategy $\mathfrak{u}\in
\mathfrak{U}$ does not include regularity assumptions on the mapping
$\nu\mapsto\mathfrak {u}[\nu]$. This is in contrast
to~\cite{BayraktarYao11}, where a continuity condition is imposed,
enabling the authors to deal with the selection problem for strategies
in the context of a stochastic differential game and use the
traditional formulation of the value functions in terms of infima (not
essential infima) and suprema. Let us mention, however, that such
regularity assumptions may preclude existence of optimal strategies in
concrete examples; see also Remark~\ref{remuncertainVolModel}.
\end{remark}

%s3.2 #&#
\subsection{PDE for the reachability set \texorpdfstring{$\Lambda$}{Lambda}}\label{secPDEZ}

In this section, we show how the PDE for the reachability set $\Lambda$
from~(\ref{eqdefinitiondureachabilityset}) can be deduced from the
geometric dynamic programming principle of
Corollary~\ref{corGDPwithlessregularity}. This equation is stated in
terms of the indicator function of the complement of the graph
of~$\Lambda$,
\[
\chi(t,z,p):= 1 - \mathbf{1}_{\Lambda(t)} (z,p) = \cases{0, &\quad if $(z,p)
\in\Lambda(t)$,
\cr
1, &\quad\mbox{otherwise}}
\]
and its lower semicontinuous envelope
\[
\chi_*(t,z,p):=\liminf_{(t',z',p')\to(t,z,p)} \chi\bigl(t',z',p'
\bigr).
\]
Corresponding results for the case without adverse player have been
obtained in~\mbox{\cite{BouchardElieTouzi09,SonerTouzi02a}}; we extend their
arguments to account for the presence of $\nu$ and the fact that we
only have a relaxed GDP. We begin by rephrasing
Corollary~\ref{corGDPwithlessregularity} in terms of $\chi$.

%
%le3.3 #&#
\begin{lemma}\label{legdppprime}
Fix $(t,z,p)\in[0,T]\times\mathbb{R}^d\times\mathbb{R}$, and let
$O\subset[0,T]\times\mathbb{R}^d\times\mathbb{R}$ be
a bounded open set containing $(t,z,p)$.
\begin{longlist}[(GDP2$_{\chi}$)]
\item[(GDP1$_{\chi}$)] Assume that $\chi(t,z,p+\varepsilon)=0$ for
    some $\varepsilon>0$. Then there
    exist $\mathfrak{u}\in\mathfrak{U}$ and $\{\alpha^{\nu }, \nu\in
    \mathcal{V}\} \subset\mathcal{A}$ such that
\[
\chi_* \bigl(\tau^\nu,Z^{\mathfrak{u},\nu}_{t,z}\bigl(
\tau^\nu \bigr),P^{\alpha^{\nu
}}_{t,p}\bigl(
\tau^\nu\bigr) \bigr)=0 \qquad\mathbb{P}\mbox{-a.s. for all $\nu\in
\mathcal{V}$,}
\]
where $\tau^{\nu}$ denotes the first exit time of $ (\cdot,
Z^{\mathfrak{u},\nu}_{t,z}, P^{\alpha^\nu}_{t,p} )$ from $O$.

\item[(GDP2$_{\chi}$)] Let $\varphi$ be a continuous
    function such that $\varphi\ge\chi$ and let
    $(\mathfrak{u},\mathfrak{a})\in\mathfrak{U}\times\mathfrak {A}$ and
    $\eta>0$ be such that
%
%
%e3.3 #&#
\begin{equation}
\label{eqgdp2chi} \varphi \bigl( \tau^{\nu}, Z^{\mathfrak{u},\nu}_{t,z}
\bigl(\tau^{\nu
}\bigr),P^{\mathfrak{a}[\nu
]}_{t,p}\bigl(
\tau^{\nu}\bigr) \bigr) \le1-\eta\qquad\mathbb{P}\mbox {-a.s. for
all $\nu \in\mathcal{V}$,}
\end{equation}
where $\tau^{\nu}$ denotes the first exit time of $ (\cdot,
Z^{\mathfrak{u},\nu}_{t,z}, P^{\mathfrak{a}[\nu]}_{t,p} )$ from $O$.
Then $\chi (t,z,\allowbreak p-\varepsilon )=0$ for all $\varepsilon>0$.
\end{longlist}
\end{lemma}

\begin{pf}
After observing that $(z,p+\varepsilon)\in\Lambda(t)$ if and only if
$\chi (t,z,p+\varepsilon)=0$ and that $(z,p)\in\bar\Lambda(t)$ implies
$\chi _{*}(t,z,p)=0$, (GDP1$_{\chi}$) follows from Corollary~\ref{corGDPwithlessregularity}, whose conditions are satisfied by
Proposition~\ref{propverifassumptionmarkovsde}. We now prove
(GDP2$_{\chi}$). Since $\varphi$ is continuous and $\partial O$ is
compact, we can find $\iota>0$ such that
\[
\varphi<1\qquad\mbox{on a $\iota$-neighborhood of }\partial O \cap \{\varphi\le
1-\eta\}.
\]
As $\chi\leq\varphi$, it follows that~(\ref{eqgdp2chi}) implies
\[
\bigl(Z^{\mathfrak{u},\nu}_{t,z}\bigl(\tau^{\nu}
\bigr),M^{\nu}\bigl(\tau^{\nu
}\bigr) \bigr) \in\Lambda_{\iota}\bigl(\tau^{\nu}\bigr) \qquad\mathbb{P}\mbox{-a.s. for all } \nu\in\mathcal{V}.
\]
Now Corollary~\ref{corGDPwithlessregularity} yields that
$(z,p-\varepsilon)\in\Lambda(t)$; that is, $\chi(t,z,p-\varepsilon)=0$.
\end{pf}

Given a suitably differentiable function $\varphi=\varphi(t,z,p)$
on $[0,T]\times\mathbb{R}
^{d+1}$, we shall denote by $\partial_{t} \varphi$ its derivative with
respect to $t$ and by $D\varphi$ and $D^{2}\varphi$ the Jacobian
and the
Hessian matrix with respect to $(z,p)$, respectively.
Given $u\in U$, $a\in\mathbb{R}^d$ and $v\in V$, we can then define
the Dynkin operator
\[
\mathcal{L}_{(Z,P)}^{u,a,v}\varphi:= \partial_t
\varphi+ \mu_{(Z,P)}(\cdot,u,v)^\top D \varphi+
\tfrac{1}{2}\operatorname{Tr} \bigl[\sigma_{(Z,P)}
\sigma^\top _{(Z,P)}(\cdot,u,a,v)D^2\varphi \bigr]
\]
with coefficients
\[
\mu_{(Z,P)}:=\pmatrix{\mu
\cr
0},\qquad \sigma_{(Z,P)}(\cdot,a,
\cdot):= \pmatrix{\sigma
\cr
a}.
\]
To introduce the associated relaxed Hamiltonians, we first define
the relaxed kernel
\[
\mathcal{N}_\varepsilon(z,q,v) = \bigl\{ (u,a)\in U\times\mathbb
{R}^d\dvtx  \bigl\llvert \sigma^\top _{(Z,P)}(z,u,a,v)q
\bigr\rrvert \le\varepsilon \bigr\},\qquad\varepsilon\ge0
\]
for $z\in\mathbb{R}^d$, $q\in\mathbb{R}^{d+1}$ and $v\in V$, as well as
the set $N_{\mathrm{Lip}}(z,q)$ of all continuous functions
\[
(\hat u,\hat a)\dvtx  \mathbb{R}^d\times\mathbb{R}^{d+1}\times V
\to U\times\mathbb{R}^{d},\qquad\bigl(z',q',v'
\bigr) \mapsto(\hat u,\hat a) \bigl(z',q',v'
\bigr)
\]
that are locally Lipschitz continuous in $(z',q')$, uniformly in $v'$
and satisfy
\[
(\hat u,\hat a)\in\mathcal{N}_{0}\qquad\mbox{on } B\times V\qquad
\mbox{for some neighborhood $B$ of $(z,q)$}.
\]
The local Lipschitz continuity will be used to ensure the local
wellposedness of the SDE for a Markovian strategy defined via $(\hat
u,\hat a)$. Setting
\[
F(\Theta,u,a,v):= \bigl\{ -\mu_{(Z,P)}(z,u,v)^\top q -
\tfrac{1}{2}\operatorname{Tr} \bigl[\sigma_{(Z,P)}\sigma
^\top _{(Z,P)}(z,u,a,v)A \bigr] \bigr\}
\]
for $\Theta=(z,q,A)\in\mathbb{R}^d\times\mathbb{R}^{d+1}\times
\mathbb{S}^{d+1}$ and $(u,a,v)\in U\times\mathbb{R}
^{d}\times V$, we can then define the relaxed Hamiltonians
%
%
%e3.4 #&#
%e3.5 #&#
\begin{eqnarray}
H^* (\Theta) &:=& \inf_{v\in V} \limsup_{\varepsilon\searrow
0,\Theta
'\rightarrow\Theta}\ \sup_{(u,a)\in\mathcal{N}_{\varepsilon
}(\Theta',v)} F\bigl(\Theta',u,a,v\bigr),
\label{eqHupper}
\\
H_* (\Theta) &:=& \sup_{(\hat u,\hat a)\in N_{\mathrm{Lip}}(\Theta)} \inf_{v\in
V} F
\bigl(\Theta,\hat u(\Theta,v),\hat a(\Theta,v),v\bigr)\label{eqHlower}.
\end{eqnarray}
[In~(\ref{eqHlower}), it is not necessary to take the relaxation
$\Theta'\rightarrow\Theta$ because $\inf_{v\in V} F$ is already lower
semicontinuous.] The question whether $H^*=H_*$ is postponed to the
monotone setting of the next section; see
Remark~\ref{remloclipselectors}.

We are now in the position to derive the PDE for $\chi$; in the
following, we write $H^*\varphi(t,z,p)$ for $H^*(z,D\varphi
(t,z,p),D^{2}\varphi
(t,z,p))$, and similarly for~$H_*$.

%
%th3.4 #&#
\begin{theorem}\label{thmpdederivavecz}
The function $ \chi_*$ is a viscosity supersolution on $[0,T)\times
\mathbb{R}^{d+1}$~of
\[
\bigl(-\partial_{t} + H^*\bigr) \varphi\ge0.
\]
The function $\chi^*$ is a viscosity subsolution on $[0,T)\times
\mathbb{R}^{d+1}$ of
\[
(-\partial_{t} + H_* )\varphi\le0.
\]
\end{theorem}

\begin{pf}
\textit{Step} 1. $\chi_*$ is a viscosity supersolution.

Let $(t_o,z_o,p_o) \in[0,T)\times\mathbb{R}^d\times\mathbb{R}$,
and let $\varphi$ be a smooth
function such that
%
%
%e3.6 #&#
\begin{equation}
\label{eqminlocstrict} (\mbox{strict}) \min_{[0,T)\times\mathbb{R}^d\times\mathbb{R}} ( \chi_* - \varphi )
= ( \chi_* - \varphi ) (t_o,z_o,p_o) = 0.
\end{equation}
We suppose that
%
%
%e3.7 #&#
\begin{equation}
\label{eqsupersolAssumpForConstra} \bigl(-\partial_{t}+H^*\bigr)\varphi(t_o,z_o,p_o)
\le-2 \eta<0
\end{equation}
for some $\eta> 0$ and work toward a contradiction.
Using the continuity of $\mu$ and $\sigma$
and the definition of the upper-semicontinuous operator $H^*$,
we can find $v_o\in V$ and $\varepsilon> 0$ such that
%
%e3.8 #&#
%e3.9 #&#
\begin{eqnarray}
\hspace*{-60pt}-\mathcal{L}_{(Z,P)}^{u,a,v_o}
\varphi(t,z,p) \le-\eta
\nonumber\\[-8pt]\label{eqpdederiv1}\\[-8pt]
\eqntext{\mbox{for all } (u,a) \in\mathcal{N}_{\varepsilon} \bigl( z,D
\varphi(t,z,p),v_o \bigr)\mbox{ and }(t,z,p) \in B_{\varepsilon},}
\end{eqnarray}
where $B_\varepsilon:=B_\varepsilon(t_o,z_o,p_o)$ denotes the open ball
of radius $\varepsilon$ around $(t_o,z_o,p_o)$. Let
\[
\partial B_\varepsilon:= \{t_o+\varepsilon\}\times\overline
{B_\varepsilon(z_o,p_o)} \cup [t_o,t_o+
\varepsilon)\times\partial B_\varepsilon(z_o,p_o)
\]
denote the parabolic boundary of $B_\varepsilon$, and set
\[
\zeta:= \min_{\partial B_\varepsilon} (\chi_* - \varphi).
\]
In view of (\ref{eqminlocstrict}), we have $\zeta>0$.

Next, we claim that there exists a sequence $(t_n,z_n,p_n,\varepsilon
_{n})_{n\ge1}$ $\subset$ $ B_\varepsilon\times(0,1)$ such that
%
%e3.10 #&#
\begin{eqnarray}\label{eqchi*entnegal0}
&&\hspace*{10pt} (t_n,z_n,p_n, \varepsilon_{n})\to(t_o,z_o,p_o,0)
\quad\mbox{and}
\nonumber\\[-8pt]\\[-8pt]
&&\hspace*{-5pt}\chi (t_n,z_n,p_n+\varepsilon_{n})=0\qquad\mbox{for all }n\geq1.\nonumber
\end{eqnarray}
In view of $\chi\in\{0,1\}$, it suffices to show that
%
%
%e3.11 #&#
\begin{equation}
\label{eqchiStarZero} \chi_*(t_o,z_o,p_o)=0.
\end{equation}
Suppose that $\chi_*(t_o,z_o,p_o)>0$; then the lower semicontinuity of
$\chi_*$ yields that $\chi_* >0$ and therefore $\chi=1$ on a
neighborhood of $(t_o,z_o,p_o)$, which implies that $\varphi$ has a strict
local maximum in $(t_o,z_o,p_o)$ and thus
\[
\partial_t \varphi(t_o,z_o,p_o)
\le0, \qquad D \varphi(t_o,z_o,p_o) = 0,\qquad
D^{2} \varphi(t_o,z_o,p_o)\le0.
\]
This clearly contradicts (\ref{eqpdederiv1}), and so the claim follows.

For any $n\geq1$, the equality in (\ref{eqchi*entnegal0}) and
(GDP1$_{\chi}$) of Lemma~\ref{legdppprime} yield
$\mathfrak{u}^n\in\mathfrak{U}$ and $\{\alpha^{n,\nu},\nu\in \mathcal
{V}\} \subset\mathcal{A}$ such that
%
%e3.12 #&#
\begin{equation}
\label{eqZetPenthetadansV} \chi_* \bigl( t\wedge\tau_n,Z^n(t\wedge
\tau_n),P^n(t\wedge\tau _n) \bigr)=0, \qquad t
\ge t_n,
\end{equation}
where
\[
\bigl(Z^n(s),P^n(s) \bigr):= \bigl(
Z^{\mathfrak{u}
^n,v_{o}}_{t_n,z_n}(s),P^{\alpha^{n,v_o}}_{t_n,p_n}(s) \bigr)
\]
and
\[
\tau_n:= \inf \bigl\{ s\ge t_n\dvtx  \bigl(s,Z^n(s),P^n(s)
\bigr) \notin B_\varepsilon \bigr\}.
\]
(In the above, $v_{o}\in V$ is viewed as a constant element of
$\mathcal{V}$.)
By (\ref{eqZetPenthetadansV}), (\ref{eqminlocstrict}) and the
definitions of $\zeta$ and $\tau_n$,
\[
-\varphi\bigl(\cdot,Z^n,P^n \bigr) (t\wedge
\tau_n) = (\chi_*-\varphi ) \bigl(\cdot,Z^n,P^n \bigr) (t\wedge\tau_n) \ge\zeta\mathbf{1}_{ \{t\ge\tau_n \}}
\ge0.
\]
Applying It\^o's formula to $-\varphi(\cdot,Z^n,P^n )$, we deduce that
%
%e3.13 #&#
\begin{equation}
\label{eqpdederiv2} \qquad S_n(t):= S_n(0) + \int
_{t_n}^{t\wedge\tau_n} \delta_n(r) \,dr + \int
_{t_n}^{t\wedge\tau_n} \Sigma_n(r)
\,dW_r \ge- \zeta\mathbf{1}_{ \{t<\tau_n \}},
\end{equation}
where
\begin{eqnarray*}
S_n(0) &:=& -\zeta-\varphi(t_n,z_n,p_n),
\\
\delta_n(r) &:=& -\mathcal{L}_{(Z,P)}^{\mathfrak{u}^n_r[v_o],\alpha
^{n,v_o}_{r},v_{o}}
\varphi \bigl(r,Z^n(r),P^n(r) \bigr),
\\
\Sigma_n(r) &:=& - D\varphi \bigl(r,Z^n(r),P^n(r)
\bigr)^\top \sigma _{(Z,P)} \bigl(Z^n(r),
\mathfrak{u}^n_r[v_{o}], \alpha
^{n,v_o}_{r},v_{o} \bigr).
\end{eqnarray*}
Define the set
\[
A_n:= [\![t_n,\tau_n]\!] \cap\{
\delta_n >-\eta\};
\]
then (\ref{eqpdederiv1}) and the definition of $\mathcal
{N}_\varepsilon$ imply that
%
%
%e3.14 #&#
\begin{equation}
\label{eqSetAnIneq} |\Sigma_n| >\varepsilon\qquad\mbox{on }A_n.
\end{equation}

%
%le3.5 #&#
\begin{lemma}\label{letrueMart}
After diminishing $\varepsilon>0$ if necessary, the stochastic exponential
\[
E_n(\cdot) = \mathcal{E} \biggl( -\int_{t_n}^{\cdot\wedge\tau_n}
\frac{\delta
_n(r)}{|\Sigma_n(r)|^2} \Sigma_n(r) \mathbf{1}_{A_n}(r)
\,dW_r \biggr)
\]
is well defined and a true martingale for all $n\geq1$.
\end{lemma}

This lemma is proved below; it fills a gap in the previous literature.
Admitting its result for the moment, integration by parts yields
\begin{eqnarray*}
(E_n S_n) (t\wedge\tau_n) &=&
S_n(0) + \int_{t_n}^{t\wedge\tau_n}
E_n \delta_n \mathbf{1}_{A_n^{c}} \,dr
\\
&&{} + \int_{t_n}^{t\wedge\tau_n} E_n \biggl(
\Sigma_n - S_n \frac{\delta_n}{|\Sigma_n|^2}\Sigma_n
\mathbf{1}_{A_n} \biggr) \,dW.
\end{eqnarray*}
As $E_n\geq0$, it then follows from the definition of $A_{n}$ that $E_n
\delta_n \mathbf{1}_{A_n^{c}}\leq0$ and so $E_n S_n$ is a local
supermartingale; in fact, it is a true supermartingale since it is
bounded from below by the martingale $-\zeta E_n$. In view
of~(\ref{eqpdederiv2}), we deduce that
\[
-\zeta-\varphi(t_n,z_n,p_n)
=(E_nS_n) (t_n)\ge\mathbb{E}
\bigl[(E_nS_n) (\tau_n) \bigr] \ge -\zeta
\mathbb{E} \bigl[\mathbf{1}_{ \{\tau_n<\tau_n \}} E_n(\tau_n)
\bigr]=0,
\]
which yields a contradiction due to $\zeta>0$ and the fact that,
by~(\ref{eqchiStarZero}),
\[
\varphi(t_n,z_n,p_n)\to\varphi(t_{o},z_{o},p_{o})=
\chi _*(t_{o},z_{o},p_{o})=0.
\]

\textit{Step} 2. $\chi^*$ is a viscosity subsolution.%\vspace*{6pt}

Let $(t_o,z_o,p_o)\in[0,T)\times\mathbb{R}^d\times\mathbb{R}$
and let $\varphi$ be a smooth
function such that
\[
\max_{[0,T)\times\mathbb{R}^d\times\mathbb{R}} \bigl(\chi^*-\varphi \bigr) = \bigl(
\chi^*-\varphi\bigr) (t_o,z_o,p_o) = 0.
\]
In order to prove that $(-\partial_{t}+H_*)\varphi(t_o,z_o,p_o)\le0$,
we assume for contradiction that
%
%e3.15 #&#
\begin{equation}
\label{eqcontrasoussol} (-\partial_{t}+H_*)\varphi(t_o,z_o,p_o)>0.
\end{equation}
An argument analogous to the proof of~(\ref{eqchiStarZero}) shows that
$\chi^*(t_o,z_o,p_o)=1$. Consider a sequence
$(t_n,z_n,p_n,\varepsilon_{n})_{n\ge1}$ in
$[0,T)\times\mathbb{R}^d\times \mathbb{R}\times(0,1)$ such that
\[
(t_n,z_n,p_n-\varepsilon_{n},
\varepsilon_{n}) \to(t_o,z_o,p_o,0)
\]
and
\[
\chi(t_n,z_n,p_n-
\varepsilon_{n}) \to\chi^* (t_o,z_o,p_o)=1.
\]
Since $\chi$ takes values in $\{0,1\}$, we must have
%
%e3.16 #&#
\begin{equation}
\label{eqsoussoleq5} \chi(t_n,z_n,p_n-
\varepsilon_{n}) =1
\end{equation}
for all $n$ large enough. Set
\[
\tilde\varphi(t,z,p):=\varphi (t,z,p)+|t-t_{o}|^{2}+|z-z_{o}|^{4}+|p-p_{o}|^{4}.
\]
Then inequality (\ref{eqcontrasoussol}) and the definition of $H_*$
imply that we can find $(\hat u,\hat a)$ in
$N_{\mathrm{Lip}}(\cdot,D\tilde\varphi)(t_o,z_o,p_o)$ such that
%
%e3.17 #&#
\begin{equation}
\label{eqsoussoleq2} %
\inf_{v\in V} \bigl( -
\mathcal{L}_{(Z,P)}^{(\hat u,\hat a)(\cdot,
D\tilde\varphi,v),v} \tilde\varphi \bigr)\ge0 \qquad\mbox{on } B_{\varepsilon}:=B_\varepsilon (t_o,z_o,p_o)
\end{equation}
for some $\varepsilon>0$. By the definition of $N_{\mathrm{Lip}}$,
after possibly changing $\varepsilon>0$, we have
%
%e3.18 #&#
\begin{equation}
\label{eqdefhatu,a} (\hat u,\hat a) (\cdot, D\tilde\varphi,\cdot) \in\mathcal
{N}_{0}(\cdot, D\tilde \varphi,\cdot) \qquad\mbox{ on }B_\varepsilon\times V.
\end{equation}
Moreover, we have
%
%e3.19 #&#
\begin{equation}
\label{eqsoussoleta} \tilde\varphi\ge\varphi+\eta\qquad\mbox{on }\partial
B_{\varepsilon}
\end{equation}
for some $\eta>0$. Since $\tilde\varphi(t_n,z_n,p_n)\to\varphi
(t_o,z_o,p_o)=\chi^*(t_o,z_o,p_o)=1$,
we can find $n$ such that
%
%e3.20 #&#
\begin{equation}
\label{eqsoussoleq6} \tilde\varphi(t_n,z_n,p_n)
\leq1+\eta/2
\end{equation}
and such that (\ref{eqsoussoleq5}) is satisfied. We fix this $n$ for
the remainder of the proof.

For brevity, we write $(\hat u,\hat a)(t,z,p,v)$ for $(\hat u,\hat
a)(z,D\tilde\varphi(t,z,p),v)$ in the sequel. Exploiting the definition
of $N_{\mathrm{Lip}}$, we can then define the mapping
$(\hat\mathfrak{u},\hat\mathfrak{a})[\cdot]\dvtx \mathcal{V}\to
\mathcal{U}\times\mathcal{A}$ implicitly via
\[
(\hat\mathfrak{u}, \hat\mathfrak{a})[\nu] = (\hat u, \hat a) \bigl(
\cdot,Z^{\hat\mathfrak{u}[\nu],\nu
}_{t_{n},z_{n}},P^{\hat\mathfrak{a}[\nu]}_{t_{n},p_{n}},\nu
\bigr)\mathbf {1}_{[t_{n},\tau^{\nu}]},
\]
where
\[
\tau^{\nu}:= \inf \bigl\{ r\ge t_n\dvtx  \bigl(
r,Z^{\hat\mathfrak
{u}[\nu],\nu
}_{t_{n},z_{n}}(r),P^{\hat\mathfrak{a}[\nu]}_{t_{n},p_{n}} (r)
\bigr) \notin B_\varepsilon \bigr\}.
\]
We observe that $\hat\mathfrak{u}$ and $\hat\mathfrak{a}$ are
nonanticipating; that is,
$(\hat\mathfrak{u},\hat\mathfrak{a})\in\mathfrak{U}\times
\mathfrak{A}$. Let us write $(Z^{\nu},P^{\nu})$ for
$(Z^{\hat\mathfrak{u},\nu
}_{t_{n},z_{n}},P^{\hat\mathfrak{a}[\nu]}_{t_{n},p_{n}})$ to alleviate
the notation. Since $\chi\leq\chi^*\leq\varphi$, the continuity of the
paths of $Z^\nu $ and $P^\nu$ and~(\ref{eqsoussoleta}) lead to
\[
\varphi \bigl( \tau^{\nu},Z^{\nu}\bigl(
\tau^{\nu}\bigr),P^{\nu}\bigl(\tau ^{\nu}\bigr) \bigr)
\le\tilde\varphi \bigl( \tau^{\nu},Z^{\nu}\bigl(
\tau^{\nu}\bigr),P^{\nu
}\bigl(\tau^{\nu}\bigr) \bigr) -
\eta.
\]
On the other hand, in view of~(\ref{eqsoussoleq2})
and~(\ref{eqdefhatu,a}), It\^o's formula applied to $\tilde\varphi$ on
$[t_n,\tau ^\nu]$ yields that
\[
\tilde\varphi \bigl( \tau^{\nu},Z^{\nu}\bigl(
\tau^{\nu}\bigr),P^{\nu
}\bigl(\tau^{\nu}\bigr) \bigr)
\le\tilde\varphi ( t_n,z_n,p_n ).
\]
Therefore, the previous inequality and (\ref{eqsoussoleq6}) show that
\[
\varphi \bigl( \tau^{\nu},Z^{\nu}\bigl(\tau^{\nu}
\bigr),P^{\nu}\bigl(\tau ^{\nu}\bigr) \bigr) \le\tilde\varphi (
t_n,z_n,p_n ) -\eta\le1-\eta/2.
\]
By (GDP2$_{\chi}$) of Lemma~\ref{legdppprime}, we deduce that $\chi
(t_{n},z_{n},p_{n}-\varepsilon_{n})=0$, which
contradicts~(\ref{eqsoussoleq5}).
\end{pf}

To complete the proof of the theorem, we still need to show
Lemma~\ref{letrueMart}. To this end, we first make the following
observation.

%
%le3.6 #&#
\begin{lemma}\label{leBMO}
Let $\alpha\in L^2_{\mathrm{loc}}(W)$ be such that $M=\int\alpha \,dW$ is a
bounded martingale and let $\beta$ be an $\mathbb{R}^d$-valued,
progressively measurable process such that $|\beta|\leq c(1+|\alpha|)$
for some constant~$c$. Then the stochastic exponential
$\mathcal{E}(\int\beta  \,dW)$ is a true martingale.
\end{lemma}

\begin{pf}
The assumption clearly implies that $\int_0^T |\beta_s|^2 \,ds<\infty$
$\mathbb{P}$-a.s. Since $M$ is bounded, we have in particular that
$M\in \mathit{BMO}$; that is,
\[
\sup_{\tau\in\mathcal{T}_0} \biggl\llVert \mathbb{E} \biggl[ \int
_\tau ^T |\alpha_s|^2 \,ds
\Big|\mathcal{F} _\tau \biggr]\biggr\rrVert _\infty<\infty.
\]
In view of the assumption, the same holds with $\alpha$ replaced by
$\beta$, so that $\int\beta \,dW$ is in $\mathit{BMO}$. This implies that
$\mathcal{E} (\int\beta \,dW)$ is a true martingale; cf.
\cite{Kazamaki94}, Theorem~2.3.
\end{pf}

\begin{pf*}{Proof of Lemma~\ref{letrueMart}}
Consider the process
\[
\beta_n(r):=\frac{\delta_n(r)}{|\Sigma_n(r)|^2} \Sigma_n(r)
\mathbf{1}_{A_n}(r);
\]
we show that
%
%
%e3.21 #&#
\begin{equation}
\label{eqbetaClaim} |\beta_n|\leq c\bigl(1+\bigl|\alpha^{n,v_o}\bigr|\bigr)
\qquad\mbox{on } [\![t_n,\tau _n]\!]
\end{equation}
for some $c>0$. Then the result will follow by applying
Lemma~\ref{leBMO} to $\alpha^{n,v_o}\mathbf{1}_{[\![t_n,\tau_n]\!]}$;
note that the stochastic integral of this process is bounded by the
definition of $\tau_n$. To prove~(\ref{eqbetaClaim}), we distinguish
two cases.

\textit{Case} 1.
$\partial_p\varphi(t_o,z_o,p_o)\neq0$. Using that $\mu$ and $\sigma$
are continuous and that $U$~and~$B_\varepsilon$ are bounded, tracing
the definitions yields that
\[
|\delta_n|\leq c \bigl\{1+ \bigl|\alpha^{n,v_o}\bigr| + \bigl|
\alpha^{n,v_o}\bigr|^2 \bigl|\partial _{pp}\varphi\bigl(
\cdot,Z^n,P^n\bigr)\bigr| \bigr\}\qquad\mbox{on } [\![t_n,\tau _n]\!],
\]
while
\[
|\Sigma_n|\geq-c + \bigl|\alpha^{n,v_o}\bigr|\bigl|\partial_{p}
\varphi\bigl(\cdot,Z^n,P^n\bigr)\bigr|\qquad\mbox{on } [\![t_n,\tau_n]\!]
\]
for some $c>0$. Since $\partial_p\varphi(t_o,z_o,p_o)\neq0$ by
assumption, $\partial_p\varphi$ is uniformly bounded away from zero on
$B_\varepsilon$, after diminishing $\varepsilon
>0$ if necessary. Hence, recalling~(\ref{eqSetAnIneq}), there is a
cancelation between $|\delta_n|$ and $|\Sigma_n|$ which allows us to
conclude~(\ref{eqbetaClaim}).

\textit{Case} 2.
$\partial_p\varphi(t_o,z_o,p_o)=0$. We first observe that
\[
\delta_n^+\leq c\bigl(1+ \bigl|\alpha^{n,v_o}\bigr|\bigr) -
c^{-1}\bigl|\alpha^{n,v_o}\bigr|^2 \partial_{pp}
\varphi\bigl(\cdot,Z^n,P^n\bigr) \qquad\mbox{on } [\![t_n,\tau _n]\!]
\]
for some $c>0$. Since $\delta_n^-$ and $|\Sigma_n|^{-1}$ are uniformly
bounded on $A_n$, it therefore suffices to show that $\partial
_{pp}\varphi \geq0$ on $B_\varepsilon$. To see this, we note
that~(\ref{eqsupersolAssumpForConstra}) and the relaxation in the
definition~(\ref{eqHupper}) of $H^*$ imply that there exists $\iota>0$
such that, for some $v\in V$ and all small $\varepsilon>0$,
%
%
%e3.22 #&#
\begin{equation}
\label{eqProofTrueMartCase2}
\qquad\qquad -\partial_t\varphi(t_o,z_o,p_o)
+ F\bigl(\Theta^\iota, u,a,v\bigr)\leq-\eta \qquad\mbox{for all $(u,a)
\in\mathcal{N}_\varepsilon\bigl(\Theta^\iota\bigr)$},
\end{equation}
where $\Theta^\iota= (z_0,p_0, D\varphi,A^\iota)$ and $A^\iota$ is the
same matrix as $D^2\varphi(t_o,z_o,p_o)$ except that the entry
$\partial _{pp}\varphi(t_o,z_o,p_o)$ is replaced by
$\partial_{pp}\varphi (t_o,z_o,p_o)-\iota$. Going back to the
definition of $\mathcal{N}_\varepsilon$, we observe that
$\mathcal{N}_\varepsilon (\Theta^\iota)$ does not depend on $\iota$
and, which is the crucial part, the assumption that
$\partial_p\varphi(t_o,z_o,p_o)=0$ implies that
$\mathcal{N}_\varepsilon(\Theta^\iota)$ is of~the form $\mathcal
{N}^U\times\mathbb{R}^d$; that is, the variable $a$ is unconstrained.
Now~(\ref{eqProofTrueMartCase2}) and the last observation show that
\[
-\bigl(\partial_{pp}\varphi(t_o,z_o,p_o)-
\iota\bigr)|a|^2 \leq c\bigl(1+|a|\bigr)
\]
for \textit{all} $a\in\mathbb{R}^d$, so we deduce that $\partial
_{pp}\varphi (t_o,z_o,p_o)\geq\iota>0$. Thus, after diminishing
$\varepsilon>0$ if necessary, we have $\partial_{pp}\varphi\geq0$ on
$B_\varepsilon$ as desired. This completes the proof.
\end{pf*}

%s3.3 #&#
\subsection{PDE in the monotone case}\label{secPDEx,y}

We now specialize the setup of Section~\ref{secBMcasemainassumptions}
to the case where the state process $Z$ consists of a pair of processes
$(X,Y)$ with values in $\mathbb{R}^{d-1}\times\mathbb{R}$, and the loss
function
\[
\ell\dvtx  \mathbb{R}^{d-1}\times\mathbb{R}\to\mathbb{R},\qquad (x,y)\mapsto
\ell(x,y)
\]
is nondecreasing in the scalar variable $y$. This setting, which was
previously studied in~\cite{BouchardElieTouzi09} for the case without
adverse control, will allow for a more explicit description of $\Lambda
$ which is particularly suitable for applications in mathematical
\mbox{finance}.

For $(t,x,y)\in[0,T]\times\mathbb{R}^{d-1}\times\mathbb{R}$ and
$(\mathfrak{u},\nu)\in\mathcal{U}\times\mathcal{V}$, let
$Z^{\mathfrak{u},\nu}_{t,x,y}=(X^{\mathfrak{u},\nu
}_{t,x},Y^{\mathfrak{u},\nu}_{t,x,y})$ be the strong solution of
(\ref{eqedsz}) with
\[
\mu(x,y,u,v):=\pmatrix{\mu_{X}(x,u,v)
\cr
\mu_{Y}(x,y,u,v)},
\qquad \sigma(x,y,u,v):=\pmatrix{\sigma_{X}(x,u,v)
\cr
\sigma_{Y}(x,y,u,v)},
\]
where $\mu_{Y}$ and $\sigma_{Y}$ take values in $\mathbb{R}$ and
$\mathbb{R}^{1\times d}$, respectively. The assumptions from
Section~\ref{secBMcasemainassumptions} remain in force; in particular,
the continuity and growth assumptions on $\mu$ and $\sigma$. In this
setup, we can consider the real-valued function
\[
{\gamma}(t,x,p):=\inf\bigl\{y\in\mathbb{R}\dvtx  (x,y,p)\in\Lambda(t)\bigr\}.
\]
In mathematical finance, this may describe the minimal capital $y$ such
that the given target can be reached by trading in the securities
market modeled by $X^{\mathfrak{u},\nu}_{t,x}$; an illustration is
given in the subsequent section. In the present context,
Corollary~\ref{corGDPwithlessregularity} reads as follows.

%
%le3.7 #&#
\begin{lemma}\label{leGDDvarpi}
Fix $(t,x,y,p)\in[0,T]\times\mathbb{R}^{d-1}\times\mathbb
{R}\times\mathbb{R}$, let $O\subset[0,T]\times\mathbb{R}
^{d-1}\times\mathbb{R}\times\mathbb{R}$ be a bounded open set
containing $(t,x,y,p)$ and
assume that $\gamma$ is locally bounded.
\begin{longlist}[(GDP2$_{\gamma}$)]
\item[(GDP1$_{\gamma}$)]
Assume that
    $y>\gamma(t,x,p+\varepsilon)$ for some $\varepsilon>0$. Then there
    exist $\mathfrak{u}\in\mathfrak{U}$ and $\{\alpha^{\nu }, \nu\in
    \mathcal{V}\} \subset\mathcal{A}$ such that
\[
Y^{\mathfrak{u},\nu}_{t,x,y}\bigl(\tau^\nu\bigr)\ge
\gamma_{*} \bigl(\tau,X^{\mathfrak{u},\nu
}_{t,x}\bigl(
\tau^\nu\bigr),P^{\alpha^{\nu}}_{t,p}\bigl(
\tau^\nu\bigr) \bigr) \qquad \mathbb{P}\mbox{-a.s. for all }
\nu\in\mathcal{V},
\]
where $\tau^{\nu}$ is the first exit time of $(\cdot, X^{\mathfrak
{u},\nu
}_{t,x}, Y^{\mathfrak{u},\nu}_{t,x,y},P^{\alpha^{\nu}}_{t,p})$ from $O$.

\item[(GDP2$_{\gamma}$)] Let $\varphi$ be a
    continuous function such that $\varphi\ge\gamma $ and let
    $(\mathfrak{u},\mathfrak{a})\in\mathfrak{U}\times\mathfrak{A}$ and
    $\eta>0$ be such that
\[
Y^{\mathfrak{u},\nu}_{t,x,y}\bigl(\tau^{\nu}\bigr)\ge\varphi \bigl(
\tau ^{\nu},X^{\mathfrak{u},\nu
}_{t,x}\bigl(\tau^{\nu}
\bigr),P^{\mathfrak{a}[\nu]}_{t,p}\bigl(\tau^{\nu}\bigr) \bigr)+\eta
\qquad\mathbb{P}\mbox{-a.s. for all $\nu\in\mathcal{V}$},
\]
where $\tau^{\nu}$ is the first exit time of $(\cdot, X^{\mathfrak
{u},\nu
}_{t,x}, Y^{\mathfrak{u},\nu}_{t,x,y},P^{\mathfrak{a}[\nu]}_{t,p})$
from $O$. Then $y\ge
\gamma(t,x,p-\varepsilon)$ for all $\varepsilon>0$.
\end{longlist}
\end{lemma}

\begin{pf}
Noting that $y>\gamma(t,x,p)$ implies $(x,y,p)\in\Lambda(t)$ and
that\break
$(x,y,p)\in\Lambda(t)$ implies $y\ge\gamma(t,x,p)$, the result follows
from Corollary~\ref{corGDPwithlessregularity} by arguments similar to
the proof of Lemma~\ref{legdppprime}.
\end{pf}

The Hamiltonians $G^{*}$ and $G_{*}$ for the PDE describing $\gamma$
are defined like $H^*$ and $H_*$ in~(\ref{eqHupper})
and~(\ref{eqHlower}), but with
\begin{eqnarray*}
&& F(\Theta,u,a,v)
\\
&&\qquad:= \bigl\{ \mu_{Y}(x,y,u,v) -\mu_{(X,P)}(x,u,v)^\top
q - \tfrac
{1}{2}\operatorname{Tr} \bigl[\sigma _{(X,P)}
\sigma^\top_{(X,P)}(x,u,a,v)A \bigr] \bigr\},
\end{eqnarray*}
where $\Theta:=(x,y,q,A)\in\mathbb{R}^{d-1}\times\mathbb
{R}\times\mathbb{R}^{d}\times\mathbb{S}^d$ and
\[
\mu_{(X,P)}(x,u,a,v):=\pmatrix{\mu_{X}(x,u,v)
\cr
0}, \qquad
\sigma_{(X,P)}(x,u,a,v):=\pmatrix{\sigma_{X}(x,u,v)
\cr
a}
\]
with the relaxed kernel $\mathcal{N}_{\varepsilon}$ replaced by
\[
\mathcal{K}_{\varepsilon}(x,y,q,v):= \bigl\{ (u,a)\in U\times
\mathbb{R}^d\dvtx  \bigl\llvert \sigma _{Y}(x,y,u,v)-q^\top
\sigma_{(X,P)}(x,u,a,v)\bigr\rrvert \le\varepsilon \bigr\}
\]
and $N_{\mathrm{Lip}}$ replaced by a set $K_{\mathrm{Lip}}$, defined
like $N_{\mathrm{Lip}}$ but in terms of~$\mathcal{K}_{0}$ instead of~$\mathcal{N}_{0}$. We then have the following result for the
semicontinuous envelopes $\gamma^{*}$ and $ \gamma_{*}$ of~$\gamma$.

%
%th3.8 #&#
\begin{theorem}\label{thmpdederivavecx,y} Assume that $\gamma$
is locally bounded. Then
$ \gamma_*$ is a viscosity supersolution on $[0,T)\times\mathbb
{R}^{d-1}\times\mathbb{R}$ of
\[
\bigl(-\partial_{t} + G^*\bigr) \varphi\ge0
\]
and $\gamma^*$ is a viscosity subsolution on $[0,T)\times\mathbb
{R}^{d-1}\times\mathbb{R}$ of
\[
(-\partial_{t} + G_* )\varphi\le0.
\]
\end{theorem}

\begin{pf}
The result follows from Lemma~\ref{leGDDvarpi} by adapting the proof of
\cite{BouchardElieTouzi09}, Theorem~2.1, using the arguments from the
proof of Theorem~\ref{thmpdederivavecz} to account for the
game-theoretic setting and the relaxed formulation of the GDP. We
therefore omit the details.
\end{pf}

We shall not discuss in this generality the boundary conditions as
$t\to T$; they are somewhat complicated to state but can be deduced
similarly as in \cite{BouchardElieTouzi09}.
Obtaining a comparison theorem at the present level of generality seems
difficult, mainly due to the presence of the sets $\mathcal
{K}_{\varepsilon}$ and $K_{\mathrm{Lip}}$ (which depend on the solution
itself) and the discontinuity of the nonlinearities at
$\partial_p\varphi=0$. It seems more appropriate to treat this question
on a case-by-case basis. In fact, once $G^*=G_*$ (see also
Remark~\ref{remloclipselectors}), the challenges in proving comparison
are similar as in the case without adverse player. For that case,
comparison results have been obtained, for example,
in~\cite{BouchardVu11} for a specific setting; see also the references
therein for more examples.

%
%re3.9 #&#
\begin{remark}\label{remloclipselectors}
Let us discuss briefly the question whether $G^*=G_*$. We shall focus
on the case where $U$ is convex and the (nondecreasing) function
$\gamma $ is strictly increasing with respect to $p$; in this case, we
are interested only in test functions $\varphi$ with
$\partial_p\varphi>0$. Under this condition,
$(u,a)\in\mathcal{K}_{\varepsilon}(\cdot,\varphi,(\partial_{x}\varphi,\partial_{p}\varphi),v)$
if and only if there exists $\zeta$ with $|\zeta |\leq1$ such that $
a=(\partial_{p}\varphi)^{-1} (\sigma_{Y}(\cdot,\varphi,u,v)-\partial
_{x}\varphi^{\top} \sigma_{X}(\cdot,u,v)-\varepsilon\zeta ). $ From
this, it is not hard to see that for such functions, the relaxation
$\varepsilon\searrow0,\Theta'\rightarrow\Theta$ in~(\ref{eqHupper}) is
superfluous as the operator is already continuous, so we are left with
the question whether
\[
\inf_{v\in V} \sup_{(u,a)\in\mathcal{K}_{0}(\Theta,v)} F(\Theta,u,a,v) =
\sup_{(\hat u,\hat a)\in K_{\mathrm{Lip}}(\Theta)} \inf_{v\in V} F\bigl(\Theta,\hat u(
\Theta,v),\hat a(\Theta,v),v\bigr).
\]
The inequality ``$\geq$'' is clear. The converse inequality will hold
if, say, for each $\varepsilon>0$, there exists a locally Lipschitz
mapping $(\hat u_{\varepsilon},\hat a_{\varepsilon})\in
K_{\mathrm{Lip}}$ such that
\[
F\bigl(\cdot,(\hat u_{\varepsilon},\hat a_{\varepsilon}) (\cdot,v),v\bigr)\ge
\sup_{(u,a)\in\mathcal{K}
_{0}(\cdot,v)}F(\cdot,u,a,v) -\varepsilon\qquad\mbox{for all } v\in V.
\]
Conditions for the existence of $\varepsilon$-optimal
\textit{continuous} selectors can be found in~\cite{KuciaNowak87},
Theorem 3.2. If $(u_{\varepsilon},a_{\varepsilon})$ is an
$\varepsilon$-optimal continuous selector, the definition of
$\mathcal{K}_0$ entails that
$a_{\varepsilon}^{\top}(\Theta,v)q_{p}=-\sigma^{\top
}_{X}(x,u_{\varepsilon}
(\Theta,v),v)q_{x}+\sigma_{Y}(x,y,\break u_{\varepsilon}(\Theta,v),v)$, where
we use the notation $\Theta=(x,y,p,(q_{x}^{\top},q_{p})^{\top},A)$.
Then $u_{\varepsilon}$ can be further approximated, uniformly on
compact sets, by a locally Lipschitz function $\hat u_{\varepsilon}$.
We may restrict our attention to $q_{p}>0$; so that, if we assume that
$\sigma^{\top}$ is (jointly) locally Lipschitz, the mapping $\hat
a_{\varepsilon}^{\top}(\Theta,v):=(q_{p})^{-1} (-\sigma ^{\top
}_{X}(x,\hat u_{\varepsilon} (\Theta,v),v)q_{x}+\sigma_{Y}(x,y, \hat
u_{\varepsilon }(\Theta,v),v) )$ is locally Lipschitz, and then $(\hat
u_{\varepsilon },\hat a_{\varepsilon})$ defines a sufficiently good,
locally Lipschitz continuous selector: for all $v\in V$,
\begin{eqnarray*}
F\bigl(\cdot,(\hat u_{\varepsilon},\hat a_{\varepsilon}) (\cdot,v),v\bigr) &\ge&
F\bigl(\cdot,(u_{\varepsilon
},a_{\varepsilon}) (\cdot,v),v\bigr)-O_{\varepsilon}(1)
\\
&\ge&\sup_{(u,a)\in\mathcal{K}_{0}} F(\cdot,u,a,v) - \varepsilon-O_{\varepsilon}(1)
\end{eqnarray*}
on a neighborhood of $\Theta$, where $O_{\varepsilon}(1)\to0$ as
$\varepsilon\to0$. One can similarly discuss other cases, for example,
when $\gamma$ is strictly concave (instead of increasing) with respect
to $p$ and the mapping $(x,y,q_{x},u,v)$ $\mapsto$
$-\sigma^{\top}_{X}(x,u,v)q_{x}+\sigma _{Y}(x,y,u,v)$ is invertible in
$u$, with an inverse, that is, locally Lipschitz, uniformly in $v$.
\end{remark}

%
%s4 #&#
\section{Application to hedging under uncertainty}\label{sectpartialHedging}

In this section, we illustrate our general results in a concrete
example, and use the opportunity to show how to extend them to a case
with unbounded strategies. To this end, we shall consider a problem of
partial hedging under Knightian uncertainty. More precisely, the
uncertainty concerns the drift and volatility coefficients of the risky
asset, and we aim at controlling a function of the hedging error; the
corresponding worst-case analysis is equivalent to a game where the
adverse player chooses the coefficients. This problem is related to the
$G$-expectation from~\cite{Peng07,Peng08}, the second order target
problem from~\cite{SonerTouziZhang2010dual} and the problem of optimal
arbitrage studied in~\cite{FernholzKaratzas11}. We let
\[
V=[\underline{\mu},\overline\mu]\times[\underline\sigma, \overline\sigma]
\]
be the possible values of the coefficients, where $\underline\mu\le
0\le\overline\mu$ and $\overline\sigma\ge\underline\sigma\ge 0$.
Moreover, $U=\mathbb{R}$ will be the possible values for the investment
policy, so that, in contrast to the previous sections, $U$ is not
bounded.

The notation is the same as in the previous section, except for an
integrability condition for the strategies that will be introduced
below to account for the unboundedness of $U$; moreover, we shall
sometimes write $\nu=(\mu,\sigma)$ for an adverse control $\nu\in
\mathcal{V}$. Given $(\mu,\sigma)\in\mathcal{V}$ and
$\mathfrak{u}\in\mathfrak {U}$, the state process
$Z^{\mathfrak{u},\nu}_{t,x,y}=(X^{\nu}_{t,x},Y^{\mathfrak{u},\nu}_{t,y})$
is governed by
\[
\frac{dX_{t,x}^{\nu}(r)}{X_{t,x}^{\nu}(r)}=\mu_r \,dr +\sigma_{r}
\,dW_{r},\qquad X_{t,x}^{\nu}(t)=x
\]
and
\[
dY^{\mathfrak{u},\nu}_{t,y}(r)=\mathfrak{u}[\nu]_{r} (
\mu_r \,dr +\sigma_{r} \,dW_{r} ),\qquad Y^{\mathfrak{u},\nu}_{t,y}(t)=y.
\]
To wit, the process $X_{t,x}^{\nu}$ represents the price of a risky
asset with unknown drift and volatility coefficients $(\mu,\sigma)$,
while $Y^{\mathfrak{u},\nu}_{t,y}$ stands for the wealth process
associated to
an investment policy $\mathfrak{u}[\nu]$, denominated in monetary
amounts. (The
interest rate is zero for simplicity.) We remark that it is clearly
necessary to use strategies in this setup: even a simple stop-loss
investment policy cannot be implemented as a control.

Our loss function is of the form
\[
\ell(x,y)=\Psi\bigl(y-g(x)\bigr),
\]
where $\Psi,g\dvtx \mathbb{R}\to\mathbb{R}$ are continuous functions of
polynomial growth. The function $\Psi$ is also assumed to be strictly
increasing and concave, with an inverse \mbox{$\Psi^{-1}\dvtx
\mathbb{R}\to\mathbb{R}$}, that is, again of polynomial growth. As a
consequence, $\ell$ is continuous and~(\ref{eqassellZpoly}) is
satisfied for some $q>0$; that is,
%
%
%e4.1 #&#
\begin{equation}
\label{eqpolygrowthSec4} \bigl|\ell(z)\bigr| \le C\bigl(1+|z|^{q}\bigr),\qquad z=(x,y)\in
\mathbb{R}^2.
\end{equation}
We interpret $g(X_{t,x}^\nu(T))$ as the random payoff of a European
option written on the risky asset, for a given realization of the drift
and volatility processes, while $\Psi$ quantifies the disutility of the
hedging error $Y^{\mathfrak{u},\nu}_{t,y}(T)-g(X^\nu_{t,x}(T))$.
In this setup,
\begin{eqnarray*}
\hspace*{-5pt}&& \gamma(t,x,p)
\\
\hspace*{-5pt}&&\!\quad =\inf \bigl\{y\in\mathbb{R}\dvtx  \exists \mathfrak{u}\in\mathfrak{U}
\mbox{ s.t. } \mathbb{E} \bigl[\Psi\bigl(Y^{\mathfrak{u},\nu}_{t,y}(T)-g
\bigl(X_{t,x}^\nu (T)\bigr)\bigr)|\mathcal{F}_{t}
\bigr]\ge p\ \mathbb{P}\mbox{-a.s. } \forall \nu\in\mathcal{V} \bigr\}
\end{eqnarray*}
is the minimal price for the option allowing to find a hedging policy
such that the expected disutility of the hedging error
is controlled by $p$.

We fix a finite constant $\bar{q}> q\vee2$ and define $\mathfrak{U}$ to
be the set of mappings \mbox{$\mathfrak{u}\dvtx
\mathcal{V}\to\mathcal{U}$}
that are nonanticipating (as in Section~\ref{secexample}) and satisfy
the integrability condition
%
%e4.2 #&#
\begin{equation}
\label{eqhypnulpunifenvc} \sup_{\nu\in\mathcal{V}}\mathbb{E} \biggl[\biggl\llvert \int
_0^T \bigl|\mathfrak{u}[\nu]_r\bigr|^2
\,dr\biggr\rrvert ^{\bar{q}/2} \biggr] <\infty.
\end{equation}
The conclusions below do not depend on the choice of $\bar{q}$. The
main result of this section is an explicit expression for the price
$\gamma(t,x,p)$.

%
%th4.1 #&#
\begin{theorem}\label{thmexplicitexpression}
Let $(t,x,p)\in[0,T]\times(0,\infty)\times\mathbb{R}$. Then
$\gamma(t,x,p)$ is finite
and given by
%
%
%e4.3 #&#
%e4.4 #&#
\begin{eqnarray}
\gamma(t,x,p)= \sup_{\nu\in\mathcal{V}^0 }\mathbb{E} \bigl[g
\bigl( X^\nu_{t,x}(T) \bigr) \bigr]+\Psi^{-1}(p)
\nonumber\\[-12pt]\label{eqsuperhedDuality} \\[-8pt]
\eqntext{\mbox{where } \mathcal{V}^0=\bigl\{(\mu,\sigma)\in\mathcal{V}\dvtx
\mu \equiv0\bigr\}.}
\end{eqnarray}
\end{theorem}

In particular, $\gamma(t,x,p)$ coincides with the superhedging price
for the shifted option $g(\cdot)+\Psi^{-1}(p)$ in the (driftless)
uncertain volatility model for $[\underline\sigma, \overline\sigma]$;
see also below. That is, the drift uncertainty has no impact on the
price, provided that $\underline\mu\le0\le\overline\mu$. Let us
remark, in this respect, that the present setup corresponds to an
investor who knows the present and historical drift and volatility of
the underlying. It may also be interesting to study the case where only
the trajectories of the underlying (and therefore the volatility, but
not necessarily the drift) are observed. This, however, does not
correspond to the type of game studied in this paper.

%s4.1 #&#
\subsection{Proof of Theorem~\protect\texorpdfstring{\ref{thmexplicitexpression}}{4.1}}
\mbox{}
\begin{pf*}{Proof of ``$\geq$'' in~(\ref{eqsuperhedDuality})}
We may assume that $\gamma(t,x,p)<\infty$. Let $y>\gamma(t,x,p)$; then
there exists $\mathfrak{u}\in\mathfrak{U}$ such that
\[
\mathbb{E} \bigl[\Psi \bigl(Y^{\mathfrak{u},\nu}_{t,y}(T)- g \bigl(
X^\nu_{t,x}(T) \bigr) \bigr) \bigr]\ge p \qquad\mbox{for all
} \nu\in\mathcal{V}.
\]
As $\Psi$ is concave, it follows by Jensen's inequality that
\[
\Psi \bigl( \mathbb{E} \bigl[Y^{\mathfrak{u},\nu}_{t,y}(T)- g \bigl(
X^\nu_{t,x}(T) \bigr) \bigr] \bigr)\ge p \qquad\mbox{for all
} \nu\in\mathcal{V}.
\]
Since the integrability condition (\ref{eqhypnulpunifenvc}) implies
that $Y^{\mathfrak{u},\nu}_{t,y}$ is a martingale for all
\mbox{$\nu\in\mathcal{V}^0$,} we conclude that
\[
\Psi \bigl( y - \mathbb{E} \bigl[g \bigl( X^\nu_{t,x}(T)
\bigr) \bigr] \bigr)\ge p\qquad\mbox{for all } \nu\in\mathcal{V}^0
\]
and hence $y\geq\sup_{\nu\in\mathcal{V}^0 }\mathbb{E}
[g ( X^\nu_{t,x}(T)
) ]+\Psi^{-1}(p)$. As $y>\gamma(t,x,p)$ was arbitrary, the
claim follows.
\end{pf*}

We shall use Theorem \ref{thmpdederivavecx,y} to derive the missing
inequality in~(\ref{eqsuperhedDuality}). Since $U=\mathbb {R}$ is
unbounded, we introduce a sequence of approximating problems $\gamma_n$
defined like~$\gamma$, but with strategies bounded by $n$,
\[
\gamma_n (t,x,p):= \inf \bigl\{ y\in\mathbb{R}\dvtx  \exists
\mathfrak{u}\in\mathfrak{U}^n \mbox{ s.t. } \mathbb{E} \bigl[\ell
\bigl( Z^{\mathfrak{u},\nu}_{t,x,y}(T) \bigr)|\mathcal{F}_t
\bigr]\ge p\ \mathbb{P}\mbox{-a.s. } \forall \nu\in\mathcal{V} \bigr\},
\]
where
\[
\mathfrak{U}^n=\bigl\{\mathfrak{u}\in\mathfrak{U}\dvtx  \bigl|\mathfrak{u}[\nu
]\bigr|\leq n \mbox{ for all }\nu\in\mathcal{V}\bigr\}.
\]
Then clearly $\gamma_{n}$ is decreasing in $n$ and
%
%
%e4.5 #&#
\begin{equation}
\label{eqynIneq} \gamma_{n}\ge\gamma,\qquad n\geq1.
\end{equation}
%

%
%le4.2 #&#
\begin{lemma}\label{lemmaapproxcontrolbornes}
Let $(t,z)\in[0,T]\times(0,\infty)\times\mathbb{R}$, $\mathfrak
{u}\in\mathfrak{U}$, and define $\mathfrak{u}
_{n}\in\mathfrak{U}$ by
\[
\mathfrak{u}_{n}[\nu]:=\mathfrak{u}[\nu]
\mathbf{1}_{ \{\llvert \mathfrak{u}[\nu
]\rrvert \le n \}},\qquad\nu\in\mathcal{V}.
\]
Then
\[
\mathop{\operatorname{ess}\sup}_{\nu\in\mathcal{V}} \bigl\llvert \mathbb{E} \bigl[
\ell \bigl( Z^{\mathfrak{u}_{n},\nu}_{t,z}(T) \bigr)- \ell \bigl(
Z^{\mathfrak{u},\nu}_{t,z}(T) \bigr) |\mathcal{F}_{t} \bigr]
\bigr\rrvert \to 0\qquad\mbox{in $L^1$ as $n\to\infty$}.
\]
\end{lemma}

\begin{pf}
Using monotone convergence and an argument as in the proof of step~1 in
Section~\ref{secproofGDP1}, we obtain that
\begin{eqnarray*}
&& \mathbb{E} \Bigl\{\mathop{\operatorname{ess}\sup}_{\nu\in\mathcal
{V}} \bigl\llvert
\mathbb {E} \bigl[ \ell \bigl( Z^{\mathfrak{u}_{n},\nu
}_{t,z}(T) \bigr)- \ell
\bigl( Z^{\mathfrak{u},\nu}_{t,z}(T) \bigr) |\mathcal{F} _{t} \bigr]
\bigr\rrvert \Bigr\}
\\
&&\qquad =\sup_{\nu\in\mathcal{V}}\mathbb{E} \bigl\{ \bigl\llvert \ell \bigl(
Z^{\mathfrak{u}_{n},\nu}_{t,z}(T) \bigr)- \ell \bigl( Z^{\mathfrak{u},\nu}_{t,z}(T)
\bigr) \bigr\rrvert \bigr\}.
\end{eqnarray*}
Since $V$ is bounded, the Burkholder--Davis--Gundy inequalities show
that there is a universal constant $c>0$ such that
\begin{eqnarray*}
\mathbb{E} \bigl\{ \bigl\llvert Z^{\mathfrak{u}_{n},\nu}_{t,z}(T)-Z^{\mathfrak{u},\nu}_{t,z}(T)
\bigr\rrvert \bigr\} &\le& c \mathbb{E} \biggl[ \int_t^T
\bigl\llvert \mathfrak{u}[\nu ]_r-\mathfrak{u}_{n}[
\nu]_r \bigr\rrvert ^2 \,dr \biggr]^{1/2}
\\
&=& c \mathbb{E} \biggl[ \int_t^T \bigl
\llvert \mathfrak{u}[\nu]_r \mathbf{1}_{ \{|\mathfrak{u}[\nu]_r|>n \}}\bigr\rrvert
^2 \,dr \biggr]^{1/2}
\end{eqnarray*}
and hence (\ref{eqhypnulpunifenvc}) and H\"older's inequality yield
that, for any given $\delta>0$,
%
%
%e4.6 #&#
\begin{eqnarray}\label{eqapprox1}
&& \sup_{\nu\in\mathcal{V}}\mathbb{P} \bigl\{\bigl\llvert
Z^{\mathfrak{u}_{n},\nu}_{t,z}(T) - Z^{\mathfrak{u},\nu
}_{t,z}(T) \bigr
\rrvert >\delta \bigr\}
\nonumber\\[-8pt]\\[-8pt]
&&\qquad \leq\delta^{-1}\sup_{\nu\in\mathcal{V}}
\mathbb{E} \bigl\{\bigl\llvert Z^{\mathfrak{u}_{n},\nu}_{t,z}(T) -
Z^{\mathfrak{u},\nu
}_{t,z}(T) \bigr\rrvert \bigr\} \to0\nonumber
\end{eqnarray}
for $n\to\infty$. Similarly, the Burkholder--Davis--Gundy inequalities
and~(\ref{eqhypnulpunifenvc}) show that $\{|Z^{\mathfrak{u}_{n},\nu
}_{t,z}(T)|+|Z^{\mathfrak{u},\nu }_{t,z}(T)|,
\nu\in\mathcal{V},n\geq1\}$ is bounded in $L^{\bar {q}}$. This yields,
on the one hand, that
%
%
%e4.7 #&#
\begin{equation}
\label{eqapprox2} \sup_{\nu\in\mathcal{V}, n\geq1}\mathbb{P} \bigl\{\bigl\llvert
Z^{\mathfrak{u}_{n},\nu}_{t,z}(T)\bigr\rrvert + \bigl\llvert
Z^{\mathfrak{u},\nu}_{t,z}(T) \bigr\rrvert >k \bigr\}\to0
\end{equation}
for $k\to\infty$, and on the other hand, in view
of~(\ref{eqpolygrowthSec4}) and $\bar{q}>q$, that
%
%
%e4.8 #&#
\begin{equation}
\label{eqapprox3} \qquad\bigl\{\ell \bigl( Z^{\mathfrak{u}_{n},\nu}_{t,z}(T) \bigr)-
\ell \bigl( Z^{\mathfrak{u},\nu}_{t,z}(T) \bigr)\dvtx  \nu\in\mathcal{V},n\geq1
\bigr\} \qquad\mbox{is uniformly integrable.}
\end{equation}
Let $\varepsilon>0$; then (\ref{eqapprox2}) and (\ref{eqapprox3}) show
that we can choose $k>0$ such that
\[
\sup_{\nu\in\mathcal{V}}\mathbb{E} \bigl[\bigl\llvert \ell \bigl(
Z^{\mathfrak{u}_{n},\nu}_{t,z}(T) \bigr)- \ell \bigl( Z^{\mathfrak{u},\nu}_{t,z}(T)
\bigr)\bigr\rrvert \mathbf {1}_{\{|Z^{\mathfrak{u}
_{n},\nu}_{t,z}(T)| + |Z^{\mathfrak{u},\nu}_{t,z}(T) |>k\}} \bigr]<\varepsilon
\]
for all $n$. Using also that $\ell$ is uniformly continuous on $\{
|z|\leq k\}$, we thus find $\delta>0$ such that
\begin{eqnarray*}
&& \sup_{\nu\in\mathcal{V}}\mathbb{E} \bigl[\bigl\llvert \ell \bigl(
Z^{\mathfrak{u}_{n},\nu}_{t,z}(T) \bigr)- \ell \bigl( Z^{\mathfrak{u},\nu}_{t,z}(T)
\bigr)\bigr\rrvert \bigr]
\\
&&\qquad \leq2\varepsilon+ \sup_{\nu\in\mathcal{V}}\mathbb{E} \bigl[\bigl\llvert
\ell \bigl( Z^{\mathfrak{u}_{n},\nu}_{t,z}(T) \bigr)- \ell \bigl(
Z^{\mathfrak{u},\nu}_{t,z}(T) \bigr)\bigr\rrvert \mathbf{1}_{\{|
Z^{\mathfrak{u}_{n},\nu}_{t,z}(T) - Z^{\mathfrak{u},\nu}_{t,z}(T)|
>\delta\}}
\bigr].
\end{eqnarray*}
By (\ref{eqapprox1}) and (\ref{eqapprox3}), the supremum on the
right-hand side tends to zero as $n\to\infty$. This completes the proof
of Lemma~\ref{lemmaapproxcontrolbornes}.
\end{pf}

\begin{pf*}{Proof of ``$\leq$'' in~(\ref{eqsuperhedDuality})} It
follows from the polynomial growth of $g$ and the boundedness of $V$
that the right-hand side of~(\ref{eqsuperhedDuality}) is finite. Thus
the already established inequality ``$\geq$''
in~(\ref{eqsuperhedDuality}) yields that $\gamma(t,x,p)>-\infty$. We
now show the theorem under the hypothesis that $\gamma(t,x,p)<\infty$
for all $p$; we shall argue at the end of the proof that this is
automatically satisfied.%\vspace*{6pt}

\textit{Step} 1. Let $\gamma_\infty:=\inf_n \gamma _n$. Then the upper
semicontinuous envelopes of $\gamma$ and $\gamma_\infty$ coincide:
$\gamma^{*}=\gamma^{*}_{\infty}$.%\vspace*{6pt}

It follows from~(\ref{eqynIneq}) that $\gamma^{*}_{\infty}\ge
\gamma^{*}$. Let $\eta>0$ and $y> \gamma(t,x,p+\eta)$. We show that
$y\ge\gamma _n(t,x,p)$ for $n$ large; this will imply the remaining
inequality $\gamma^{*}_{\infty}\le\gamma^{*}$. Indeed, the definition
of $\gamma$ and Lemma~\ref{lemmaapproxcontrolbornes} imply that we can
find $\mathfrak{u}\in\mathfrak{U}$ and $\mathfrak{u}_{n}\in\mathfrak{U}
^n$ such that
\[
J(t,x,y, \mathfrak{u}_{n})\ge J(t,x,y,\mathfrak{u})-
\epsilon_{n}\ge p+\eta-\epsilon_{n} \qquad\mathbb{P}\mbox{-a.s.},
\]
where $\epsilon_{n}\to0$ in $L^1$. If $K_{n}$ is defined like $K$, but
with $\mathfrak{U}^{n}$ instead of $\mathfrak{U}$, then it follows that
$K_{n}(t,x,y)\ge p +\eta-\epsilon_{n}$ $\mathbb{P}$-a.s. Recalling that
$K_{n}$ is deterministic (cf.~Proposition~\ref{propverifassumptionmarkovsde}), we may replace
$\epsilon_n$ by $\mathbb{E} [\epsilon_n]$ in this inequality. Sending
\mbox{$n\to\infty$}, we then see that $\lim_{n\to\infty}K_{n}(t,x,y)\ge p
+\eta$, and therefore $K_{n}(t,x,y)\ge p +\eta/2$ for $n$ large enough.
The fact that $y\ge \gamma_n(t,x,p)$ for $n$ large then follows from
the same considerations as in
Lemma~\ref{lemselectionofepsoptimalstrat}.%\vspace*{6pt}

\textit{Step} 2. The relaxed semi-limit
\[
\bar\gamma^{*}_{\infty}(t,x,p):=\mathop{\limsup
_{n\rightarrow\infty}}_{(t',x',p')\rightarrow
(t,x,p)}\gamma^{*}_n
\bigl(t',x',p'\bigr)
\]
is a viscosity subsolution on $[0,T)\times(0,\infty)\times\mathbb
{R}$ of
%
%e4.9 #&#
\begin{equation}
\label{eqvarpisoussolbsb} -\partial_t\varphi+ \inf_{\sigma\in[\underline\sigma,\overline\sigma]}
\biggl\{ - \frac{1}2\sigma^2x^2
\partial_{xx}\varphi \biggr\} \le0
\end{equation}
and satisfies the boundary condition
$\bar\gamma^{*}_{\infty}(T,x,p)\leq g(x)+\Psi^{-1}(p)$.%\vspace*{6pt}

We first show that the boundary condition is satisfied. Fix $(x,p)\in
(0,\infty)\times\mathbb{R}$ and let $y>g(x)+\Psi^{-1}(p)$; then
$\ell(x,y)> p$. Let $(t_{n},x_n,p_n)\to(T, x,p)$ be such that $\gamma
_{n}(t_{n},x_n,p_n)\to\bar\gamma^{*}_{\infty}(T, x,p)$. We consider the
strategy $\mathfrak{u}\equiv0$ and use the arguments from the proof of
Proposition~\ref{propverifassumptionmarkovsde} to find a constant $c$
independent of $n$ such that
\[
\mathop{\operatorname{ess}\sup}_{\nu\in\mathcal{V}}\mathbb {E} \bigl[\bigl|Z^{0,\nu
}_{t_{n},x_{n},y}(T)-(x,y)\bigr|^{\bar
{q}}
|\mathcal{F}_{t_n} \bigr] \le c \bigl(|T-t_{n}|^{\bar{q}/2}+
|x-x_{n}|^{\bar{q}} \bigr).
\]
Similar to the proof of Lemma~\ref{lemmaapproxcontrolbornes}, this
implies that there exist constants $\varepsilon_{n}\to0$ such that
\[
J(t_{n},x_{n},y,0)\ge\ell(x,y)-\varepsilon_{n}
\qquad\mathbb{P}\mbox{-a.s.}
\]
In view of $\ell(x,y)>p$, this shows that $y\ge\gamma
_{n}(t_{n},x_n,p_n)$ for $n$ large enough, and hence that $y\ge\bar
\gamma^{*}_{\infty}(T, x,p)$. As a result, we have $\bar\gamma
^{*}_{\infty}(T,x,p)\leq g(x)+\Psi^{-1}(p)$.

It remains to show the subsolution property.
Let $\varphi$ be a smooth function, and let $(t_o,x_o,p_o)\in
[0,T)\times
(0,\infty)\times\mathbb{R}$ be such that
\[
\bigl(\bar\gamma^*_{\infty} - \varphi\bigr) (t_o,x_o,p_o)=
\max\bigl(\bar\gamma ^*_{\infty
}-\varphi\bigr)=0.
\]
After passing to a subsequence, \cite{Barles94}, Lemma 4.2, yields
$(t_n,x_n,p_n)\to(t_o,x_o,p_o)$ such that
\[
\lim_{n\rightarrow\infty}\bigl(\gamma_n^*-\varphi\bigr)
(t_n,x_n,p_n)=\bigl(\bar \gamma
^*_{\infty}-\varphi\bigr) (t_o,x_o,p_o)
\]
and such that $(t_n,x_n,p_n)$ is a local maximizer of $(\gamma
_n^*-\varphi )$. Applying Theorem \ref{thmpdederivavecx,y} to $\gamma
^{*}_{n}$, we deduce that
%
%e4.10 #&#
\begin{equation}
\label{eqsoussolvarpin} \sup_{(\hat u,\hat a)\in K^n_{\mathrm{Lip}}(\cdot,D\varphi)}\inf_{(\mu,\sigma)\in V} G
\varphi\bigl(\cdot,(\hat u,\hat a) (\mu,\sigma),(\mu,\sigma)\bigr)
(t_n,x_n,p_n)\le0,
\end{equation}
where
\begin{eqnarray*}
&& G\varphi\bigl(\cdot,(u,a),(\mu,\sigma)\bigr)
\\
&&\qquad:=u\mu- \partial_t
\varphi-\mu x \partial_x\varphi-\tfrac{1}2 \bigl( \sigma
^{2}x^{2} \partial_{xx}\varphi+ a^{2}
\partial_{pp}\varphi+2 \sigma x a \partial _{xp}\varphi
\bigr)
\end{eqnarray*}
and $K^{n}_{\mathrm{Lip}}(\cdot,D\varphi)(t_{n},x_{n},p_{n})$ is the
set of locally Lipschitz mappings $(\hat u,\hat a)$ with values in
$[-n,n]\times \mathbb{R}$ such that
\[
\sigma\hat u(x,q_{x},q_{p},\mu,\sigma)=x \sigma
q_{x}+q_{p} \hat a(x,q_{x},q_{p},
\mu,\sigma) \qquad\mbox{for all } \sigma\in[\underline\sigma, \overline\sigma]
\]
for all $(x,(q_{x},q_{p}))$ in a neighborhood of $(x_{n}, D\varphi
(t_{n},x_{n},p_{n}))$.
Since the mapping
\[
(0,\infty)\times\mathbb{R}^2\times[\underline\mu,\overline\mu ]
\times[\underline\sigma,\overline\sigma]\to\mathbb{R}^2
\qquad(x,q_x,q_p,\mu,\sigma )\mapsto (x q_x,0 )
\]
belongs to $K^{n}_{\mathrm{Lip}}(\cdot,D\varphi)(t_{n},x_{n},p_{n})$
for $n$ large enough, (\ref{eqsoussolvarpin}) leads to
\[
-\partial_t\varphi+ \inf_{\sigma\in[\underline\sigma,\overline
\sigma]} \biggl\{ -
\frac{1}2\sigma^2x^2\partial_{xx}
\varphi \biggr\}(t_n,x_n,p_n)\le0
\]
for $n$ large. Here the nonlinearity is continuous; therefore, sending
$n\rightarrow\infty$ yields~(\ref{eqvarpisoussolbsb}).%\vspace*{6pt}

\textit{Step} 3.
We have $\bar\gamma^{*}_{\infty} \le\pi$ on $
[0,T]\times (0,\infty)\times\mathbb{R} $, where
\[
\pi(t,x,p):= \sup_{\nu\in\mathcal{V}^0 }\mathbb{E} \bigl[g \bigl(
X^\nu_{t,x}(T) \bigr) \bigr]+\Psi^{-1}(p)
\]
is the right-hand side of~(\ref{eqsuperhedDuality}).%\vspace*{6pt}

Indeed, our assumptions on $g$ and $\Psi^{-1}$ imply that $\pi$ is
continuous with polynomial growth. It then follows by standard arguments
that $\pi$ is a viscosity supersolution on $[0,T)\times(0,\infty
)\times\mathbb{R}$ of
\[
-\partial_t\varphi+ \inf_{\sigma\in[\underline\sigma,\overline
\sigma]} \biggl\{ -
\frac{1}2\sigma^2x^2\partial_{xx}
\varphi \biggr\} \ge0
\]
and clearly the boundary condition $\pi(T,x,p)\geq g(x)+\Psi^{-1}(p)$
is satisfied. The claim then follows from step~2 by comparison.

We can now deduce the theorem: we have $\gamma\leq\gamma^{*}$ by the
definition of $\gamma^*$ and $\gamma^* = \gamma_{\infty}^{*}$ by
step~1. As $\gamma_{\infty}^{*} \le\bar\gamma_{\infty}^{*}$ by
construction, step~3 yields the result.

It remains to show that $\gamma<\infty$. Indeed, this is clearly
satisfied when $g$ is bounded from above. For the general case, we
consider $g_m=g\wedge m$ and let $\gamma_m$ be the corresponding value
function. Given $\eta>0$, we have $\gamma_m(t,x,p+\eta)<\infty$ for all
$m$ and so~(\ref{eqsuperhedDuality}) holds for $g_m$. We see
from~(\ref{eqsuperhedDuality}) that
$y:=1+\sup_m\gamma_m(t,x,\allowbreak p+\eta)$
is finite. Thus, there exist $\mathfrak{u}_m\in\mathfrak{U}$ such that
\[
\mathbb{E} \bigl[\Psi \bigl(Y^{\mathfrak{u}_m,\nu}_{t,y}(T)-
g_m \bigl( X^\nu_{t,x}(T) \bigr) \bigr)|
\mathcal{F}_t \bigr]\ge p+\eta\qquad\mbox{for all } \nu \in
\mathcal{V}.
\]
Using once more the boundedness of $V$, we see that for $m$ large enough,
\[
\mathbb{E} \bigl[\Psi \bigl(Y^{\mathfrak{u}_m,\nu}_{t,y}(T)- g \bigl(
X^\nu_{t,x}(T) \bigr) \bigr)|\mathcal{F}_t
\bigr]\ge p \qquad\mbox{for all } \nu\in \mathcal{V},
\]
which shows that $\gamma(t,x,p)\leq y<\infty$.
\end{pf*}

%
%re4.3 #&#
\begin{remark}\label{remuncertainVolModel}
We sketch a probabilistic proof for the inequality ``$\leq$'' in
Theorem~\ref{thmexplicitexpression}, for the special case without drift
($\underline\mu= \overline\mu=0$) and $\underline\sigma>0$. We focus on
$t=0$, and recall that $y_0:=\sup_{\nu\in\mathcal{V}^0 }\mathbb{E}[g(
X^\nu_{0,x}(T))]$ is the superhedging price for $g(\cdot)$ in the
uncertain volatility model. More precisely, if $B$ is the
coordinate-mapping process on $\Omega=C([0,T];\mathbb{R})$, there
exists an $\mathbb {F}^B$-progressively measurable process $\vartheta$
such that
\[
y_0+\int_0^T
\vartheta_s \frac{dB_s}{B_s} \geq g(B_T)\qquad P^\nu \mbox{-a.s. for all $\nu\in\mathcal{V}^0$,}
\]
where $P^\nu$ is the law of $X^\nu_{0,x}$ under $P$; see, for example,
\cite{NutzSoner10}. Seeing $\vartheta$ as an adapted functional of $B$,
this implies that
\[
y_0+\int_0^T
\vartheta_s\bigl(X^\nu_{0,x}\bigr)
\frac{dX^\nu
_{0,x}(s)}{X^\nu
_{0,x}(s)} \geq g\bigl(X^\nu_{0,x}(T)\bigr)\qquad P
\mbox{-a.s. for all $\nu\in \mathcal{V}^0$}.
\]
Since $X^\nu_{0,x}$ is nonanticipating with respect to $\nu$, we see
that $\mathfrak{u}[\nu]_s:=\vartheta_s(X^\nu_{0,x})$ defines a
nonanticipating strategy such that, with $y:=y_0+\Psi^{-1}(p)$,
\[
y+\int_0^T \mathfrak{u}[\nu]_{s}
\frac{dX^\nu_{0,x}(s)}{X^\nu
_{0,x}(s)} \geq g\bigl(X^\nu_{0,x}(T)\bigr) +
\Psi^{-1}(p);
\]
that is,
\[
\Psi \bigl(Y^{\mathfrak{u},\nu}_{0,y}(T)-g\bigl(X_{0,x}^\nu(T)
\bigr)\bigr)\ge p
\]
holds even $P$-almost surely, rather than only in expectation, for all
$\nu\in\mathcal{V}^0$, and $\mathcal{V}^0=\mathcal{V}$ because of our
assumption that $\underline\mu= \overline\mu=0$. In particular, we have
the \textit{existence of an optimal strategy} $\mathfrak{u}$. (We
notice that, in this respect, it is important that our definition of
strategies does not contain regularity assumptions on
$\nu\mapsto\mathfrak{u}[\nu]$.)

Heuristically, the case with drift uncertainty (i.e., $\underline\mu
\neq\overline\mu$) can be reduced to the above by a Girsanov change of
measure argument; for example, if $\mu$ is deterministic, then we can
take $\mathfrak{u}[(\mu,\sigma)]:=\mathfrak{u}[(0,\sigma^\mu)]$, where
$\sigma^\mu(\omega ):=\sigma(\omega+\int\mu_t \,dt)$. However, for
general $\mu$, there are difficulties related to the fact that a
Girsanov Brownian motion need not generate the original filtration
(see, e.g., \cite{FeldmanSmorodinsky97}), and we shall not enlarge on
this.
\end{remark}

\section*{Acknowledgments}
We are grateful to Pierre Cardaliaguet for valuable discussions and to
the anonymous referees for careful reading and helpful comments.

%suskaldyti doi

% imsref loaded by linak, 2013-11-04 16:29:30

\printaddresses

\end{document}